\documentstyle[12pt]{article}

\input amssym.def
\input amssym.tex

\def\vv{{\underline{v}}}
\def\tt{{\underline{t}}}
\def\ww{{\underline{w}}}
\def\mm{{\underline{m}}}
\def\aa{{\underline{a}}}
\def\kk{{\underline{k}}}

\def\w{\widetilde}
\def\1{{\underline{1}}}
\def\betab{\bar\beta}
\def\P{\Bbb P}
\def\R{\Bbb R}

\def\Z{\Bbb Z}

\def\C{\Bbb C}
\def\CP{\Bbb C\Bbb P}

\newtheorem{theorem}{Theorem}
\newtheorem{corollary}{Corollary}
\newtheorem{lemma}{Lemma}
\newtheorem{proposition}{Proposition}

\newenvironment{definition}
{\smallskip\noindent{\bf Definition\/}:}{\smallskip\par}
\newenvironment{statement}
{\smallskip\noindent{\bf Statement\/}.}{\smallskip\par}

\newenvironment{remark}
{\smallskip\noindent{\bf Remark\/}.}{\smallskip\par}
\newenvironment{remarks}
{\smallskip\noindent{\bf Remarks\/}.}{\smallskip\par}
\newenvironment{proof}
{\noindent{\bf Proof\/}.}{{ $\Box$}\smallskip\par}

\title{The Alexander polynomial of a plane curve singularity
via the ring of functions on it}
\author{A.Campillo
\and F.Delgado \thanks{First two authors were partially supported by
DGICYT PB97-0471 and by Junta de Castilla y Le\'on:
VA 102/01. Address:
University of Valladolid, Dept. of Algebra, Geometry and Topology,
47005 Valladolid, Spain.
E-mail: campillo\symbol{'100}agt.uva.es,
fdelgado\symbol{'100}agt.uva.es}
\and S.M.Gusein--Zade \thanks{Partially supported by grants Iberdrola,
RFBR--01--01--00739, NWO 047.008.005 and INTAS--00--0259.
Address: Moscow State University,
Dept. of Mathematics and Mechanics, Moscow, 119899, Russia.
E-mail: sabir\symbol{'100}mccme.ru}}

\date{}

\begin{document}

\def\eps{\varepsilon}

\maketitle

\begin{abstract}
We prove two formulae which express the Alexander polynomial $\Delta^C$
of several variables of a plane curve singularity $C$ in terms of the ring
${\cal O}_{C}$ of germs of analytic functions on the curve. One of them
expresses $\Delta^C$ in terms of dimensions of some factors corresponding
to a (multi-indexed) filtration on the ring ${\cal O}_{C}$. The other one
gives the coefficients of the Alexander polynomial $\Delta^C$ as Euler
characteristics of some explicitly described spaces (complements to
arrangements of projective hyperplanes).
\end{abstract}

\section{Introduction}

The ring ${\cal O}_{X}$ of germs of holomorphic functions on a germ $X$ of
an analytic set determines $X$ itself (up to analytic equivalence).
Thus all invariants of $X$, in particular, topological ones, can ``be read"
from ${\cal O}_{X}$. There arises a general problem to find expressions for
such invariants in terms of the ring ${\cal O}_{X}$.

Let $C$ be a germ of a reduced plane curve at the origin in $\C^2$ and let
$C=\bigcup\limits_{i=1}^{r}C_i$ be its representation as the union of
irreducible components (with a fixed numbering). Let
$\Delta^C(t_1, \ldots, t_r)$ be the Alexander polynomial of the link
$C\cap S_\varepsilon^3\subset S_\varepsilon^3$ for $\varepsilon > 0$ small
enough (see, e.g., {\cite{EN}}). The Alexander polynomial $\Delta^C(\tt)$
($\tt=(t_1, \ldots, t_r)$) is a complete topological invariant of a plane
curve singularity $C$ (\cite{Y}). We prove two formulae for the
Alexander
polynomial $\Delta^C$ in terms of the ring ${\cal O}_{C}$ of
germs of analytic functions on the curve $C$. For the case
of an irreducible plane curve singularity
($r=1$) the corresponding result was described in~{\cite{CDG2}}.
For the general case the result has been announced in~{\cite{CDG4}}.
A global analogue of the statement from~{\cite{CDG2}}
for plane curves with one place at infinity can be found in~{\cite{CDG3}}.

Let $\C_i$ be the complex line with the coordinate $\tau_i$ ($i=1, \ldots, r$)
and let $\varphi_i:(\C_i, 0)\to(\C^n, 0)$ be parameterizations
(uniformizations) of the branches $C_i$ of the curve $C$, i.e., germs of
analytic maps such that ${\rm{Im}}\,\varphi_i=C_i$ and $\varphi_i$ is an
isomorphism between $\C_i$ and $C_i$ outside of the origin. Let
${\cal O}_{\C^2, 0}$ be the ring of germs of holomorphic functions at the
origin in $\C^2$.
For a germ $g\in{\cal O}_{\C^2, 0}$, let $v_i=v_i(g)$ and $a_i=a_i(g)$ be
the power of the leading term and the coefficient at it in the power series
decomposition of the germ $g\circ\varphi_i:(\C_i,0)\to \C$:
$g\circ\varphi_i(\tau_i)=
a_i\cdot\tau_i^{v_i}+{~terms~of~higher~degree}$ ($a_i\ne 0$).
If $g\circ\varphi_i(t)\equiv 0$, $v_i(g)$ is assumed to be equal to $\infty$
and $a_i(g)$ is not defined. The numbers $v_i(g)$ and $a_i(g)$ are defined
for elements $g$ of the ring ${\cal O}_C$ of functions on the curve $C$
as well.

For $\vv=(v_1, \ldots, v_r) \in \Z^n$, let $J(\vv)=\{g\in {\cal O}_C:
v_i(g)\ge v_i;\, i=1, \ldots, r\}$ which is an ideal in ${\cal O}_C$,
let $c(\vv) = \mbox{dim } J(\vv)/J(\vv+{\1})$, where $\1=(1, \ldots, 1)$.
Let
$$
L_C(t_1, \ldots,
t_r) =\sum\limits_{\vv\in\Z^r}c(\vv)\cdot\tt^\vv,
$$
$$
P_C(t_1, \cdots ,t_r) =\frac{L_C(t_1, \cdots, t_r)\cdot
\prod\limits_{i=1}^r (t_i-1)}{t_1\cdot\ldots\cdot t_r-1}
$$
($\tt^\vv=t_1^{v_1}\cdot\ldots\cdot t_r^{v_r}$).
We shall show that $P_C$ is a polynomial in $t_1$, \dots, $t_r$.
Let $F_{\vv}\subset (\C^*)^r$ be the set of all $r$-tuples
$(a_1(g), \ldots, a_r(g))$ for $g\in{\cal O}_{\C^2, 0}$ with
$v_i(g)=v_i$, $i=1, \ldots, r$. The subspace $F_{\vv}$ is invariant
with respect to multiplication by non-zero complex numbers
(in fact, if $F_{\vv}$ is not empty, it is the complement to an
arrangement of hyperplanes in a vector space of dimension $c(\vv)$).
Let $\P F_{\vv}$ be the projectivization of $F_{\vv}$, i.e.,
the factor-space $F_{\vv}/\C^*$ with respect to this $\C^*$--action.

Now we formulate our main results (see more precise explanations
and definitions below).

\begin{theorem}\label{theo1}
For a plane curve singularity
$C=\bigcup\limits_{i=1}^{r}C_i\subset (\C^2, 0)$, $r>1$,
$$P_C(t_1, \ldots, t_r)=\Delta^C(t_1, \ldots, t_r).$$
\end{theorem}

\begin{theorem}\label{theo2}
For a plane curve singularity $C=\bigcup\limits_{i=1}^{r}C_i\subset
(\C^2, 0)$, $r>1$,
$$\Delta^C(t_1, \ldots, t_r)=\sum\limits_{\vv\in\Z_{\ge 0}^r}
\chi(\P F_{\vv})\cdot\tt^\vv.$$
\end{theorem}

It was found that these results can be formulated in
terms of the integral with respect to Euler characteristic over
the projectivization of the space of functions in two variables
defined in the spirit of the motivic integration; see \cite{CDG5}.
Recently W.Ebeling (\cite{Eb}) has found that there is a
relation between the Poincar\'e series of the natural filtration
in the ring of functions on a two--dimensional quasihomogeneous
hypersurface singularity and the characteristic polynomial of the
monodromy operator.
This permits to expect that the discussed connections are more
deep and must have broader field
of applications.

\section{Necessary concepts and facts.}

In this section we give more precise definitions of the objects
used in the formulation of Theorems~\ref{theo1} and \ref{theo2}
and in the proofs and describe some of their properties.

\subsection{The Alexander polynomial of an algebraic
link.}\label{s_alexander}
The Alexander polynomial (in $r$ variables) is an invariant of a link
with $r$ (numbered) components in the sphere $S^3$. The general definition
can be found, e.g., in \cite{EN}. To a plane curve singularity
$C=\bigcup\limits_{i=1}^r C_i\subset (\C^2, 0)$ there corresponds
the link $C\cap S^3_\varepsilon$ in the 3--sphere $S^3_\varepsilon$
of radius $\varepsilon$ centred at the origin in the complex plane $\C^2$
with $\varepsilon$ small enough. For such a link (an algebraic one) we rather use
not the general definition of the Alexander polynomial
$\Delta^C(t_1, \ldots, t_r)$, but a formula for
it in terms of an embedded resolution $\pi:({X}, {D})\to(\C^2, 0)$
of the curve singularity $C$.

Let the curve $C$ be given by an equation $f=0$ and let
$f=\prod\limits_{i=1}^{r}f_i$, where $f_i=0$ is an equation of
the curve $C_i$. Let
$\pi:(X, D)\to(\C^2, 0)$ be an
(embedded) resolution of the plane curve $C=\bigcup\limits_{i=1}^{r}C_i$.
Such a resolution can be described by its dual graph $\Gamma$. Vertices
of the graph $\Gamma$ correspond to components of the total transform
$(f\circ\pi)^{-1}(0)$ of the curve $C$ (i.e., to components of the
exceptional divisor $D=\pi^{-1}(0)$ of the resolution and to strict
transforms $\widetilde C_i$ of the branches $C_i$ of the curve $C$; in the
last case
they are depicted by arrows). Two vertices of the graph $\Gamma$ are
connected by an edge if the corresponding components intersect.
The graph $\Gamma$ is a tree. The starting point of the graph
(the starting divisor of the resolution) will be denoted by ${\bf 1}$.
There is a partial order on the set of vertices of the graph
$\Gamma$:  $\sigma^\prime < \sigma$ iff the geodesic in $\Gamma$ from
the vertex ${\bf 1}$ to the vertex $\sigma$ passes through the
vertex $\sigma^\prime$. A vertex $\delta$ corresponding to a
component $E_\delta$ of the exceptional divisor is said to be
{\bf a dead end} if it is connected with only one vertex (i.e., if $E_\delta$
intersects only one component of the total transform of the curve $C$).
A vertex $\sigma$ is said to be {\bf a star point} of the
resolution if it is
connected with at least three vertices. To each dead end $\delta$ in the
graph $\Gamma$ except (possibly) the vertex ${\bf 1}$ there
corresponds the nearest star point $st_\delta$ such that
$st_\delta<\delta$. All
vertices $\sigma^\prime$ such that $st_\delta<\sigma^\prime\le\delta$
form {\bf the tail} of the resolution graph corresponding to the dead
end $\delta$. A vertex $\sigma$ is said to be
{\bf a separation} {\bf point} of the graph $\Gamma$ if there
exist two branches $C_i$ and
$C_j$ of the curve $C$ such that $\sigma<\widetilde C_i$, $\sigma<
\widetilde C_j$ and $\sigma$ is the maximal vertex with these properties
(one also says that $\sigma$ is the separation point
between the branches $C_i$ and $C_j$ or between $\widetilde C_i$
and $\widetilde C_j$). Let $st_1$ be the first (i.e., the minimal) separation
point of the graph $\Gamma$. It is possible that $st_1$ is not a star point
(if $st_1={\bf 1}$). However in what follows we always include $st_1$
in the set of star vertices.

For a vertex $\sigma$ corresponding to a component $E_\sigma$ of the
exceptional
divisor (a complex projective line), let ${\stackrel{\circ}{E}}_\sigma$
be the "smooth part" of the component $E_\sigma$, i.e., $E_\sigma$
minus intersection points with other components of the total transform
of the curve $C$. These intersection points are in one-to-one
correspondence with connected components of the complement
$\mbox{$(f\circ\pi)^{-1}(0)\setminus {\stackrel{\circ}{E}}_\sigma$}$. An
intersection point is said to be {\bf essential} if the
corresponding connected
component contains a component of the strict transform of the curve $C$. Let
$s_\sigma$ be the number of essential points on the component $E_\sigma$, and
let ${\widetilde E}_\sigma$ be the complement to the set of essential points
in $E_\sigma$. Let $m_j^\sigma$ ($j=1, 2, \ldots, r$) be the multiplicity
of the lifting $f_j\circ\pi$ of the function $f_j$ (the equation of the
component $C_j$) to the space $X$ of the resolution along the component
$E_\sigma$, $\mm^\sigma:=(m_1^\sigma, \ldots, m_r^\sigma)$.

\medskip
D.Eisenbud and W.Neumann (\cite{EN}) gave a formula for the Alexander
polynomial $\Delta^C(t_1, \ldots, t_r)$ of the curve $C$ in terms of an
embedded resolution of the curve $C$.

\begin{proposition}\label{propEN} For $r>1$,
$$
\Delta^C(t_1, \ldots, t_r)=\prod\limits_{E_\sigma\subset{\cal
D}}
\left(1-\tt^{{\mm}^\sigma}\right)^{-\chi({\stackrel{\circ}{E}}_\sigma)}.\eqno(*)
$$
\end{proposition}

The formula ($*$) is an analogue of the formula of N.A'Campo for the
zeta--function of the monodromy transformation of the curve $C$.

\begin{remarks}
{\bf 1.} According to the definition, the Alexander polynomial
$\Delta^C(t_1, \ldots, t_r)$ of a link is well defined only up
to multiplication by monomials
$\pm\tt^{\underline{m}}=\pm t_1^{m_1}\cdot\ldots\cdot t_r^{m_r}$
($\tt=(t_1, \ldots, t_r)$, ${\underline{m}}=(m_1, \ldots, m_r)\in \Z^r$).
For algebraic links the formula (*) fixes the choice of the
Alexander polynomial in such a way that it is really a polynomial
(i.e., does not contain monomials with negative powers) and its value
at the origin ($\tt=0$) is equal to 1.

\smallskip\noindent {\bf 2.} There is some difference in definitions
(or rather in descriptions) of the Alexander polynomial for a curve
with one branch ($r=1$) or with many branches ($r>1$) (see, e.g.,
{\cite{EN}}). In order to have all the results (Theorems~{\ref{theo1}}
and~{\ref{theo2}}) valid for $r=1$ as well, for an irreducible
plane curve singularity $C$,
$\Delta^C(t)$ should be not the Alexander polynomial, but rather the
zeta-function of the monodromy, equal to the Alexander polynomial divided by
$(1-t)$. In this case $\Delta^C(t)$ is not a polynomial, but an infinite power
series (defined by the formula (*)). The results are valid for
this case as well. However
since the case of an irreducible plane curve singularity ($r=1$) has been
described in~\cite{CDG2}, here we shall suppose that $r>1$.
For $r>1$, $\Delta^C(t, \ldots, t)$ is nothing else but the
zeta-function of the monodromy transformation of the curve
singularity $C$.
\end{remarks}

\subsection{The extended semigroup and the Poincar\'e
polynomial of a curve singularities.}\label{s-extended}
These notions can be defined not only for plane curve singularities,
but for curve singularities in spaces of any dimension.
Let $C=\bigcup\limits_{i=1}^{r}C_i$ be a germ of a reduced curve at the
origin in $\C^n$ ($C_i$ are irreducible components (branches)
of the curve $C$). Let
$\C_i$ be the complex line with the coordinate $\tau_i$ ($i=1, \ldots, r$)
and let $\varphi_i:(\C_i, 0)\to(\C^n, 0)$ be parameterizations
(uniformizations) of the branches $C_i$ of the curve $C$, i.e., germs of
analytic maps such that ${\rm{Im}}\,\varphi_i=C_i$ and $\varphi_i$ is an
isomorphism between $\C_i$ and $C_i$ outside of the origin. Let
${\cal O}_{\C^n, 0}$ be the ring of germs of holomorphic functions at the
origin in $\C^n$.
For a germ $g\in{\cal O}_{\C^n, 0}$, let $v_i=v_i(g)$ and $a_i=a_i(g)$ be
the power of the leading term and the coefficient at it in the power series
decomposition of the germ $g\circ\varphi_i:(\C_i,0)\to \C$:
$g\circ\varphi_i(\tau_i)=
a_i\cdot\tau_i^{v_i}+{~terms~of~higher~degree}$ ($a_i\ne 0$).
If $g\circ\varphi_i(t)\equiv 0$, $v_i(g)$ is assumed to be equal to $\infty$
and $a_i(g)$ is not defined. The numbers $v_i(g)$ and $a_i(g)$ are defined
for elements $g$ of the ring ${\cal O}_C$ of functions on the curve $C$
as well.

The {\bf semigroup} $S=S_C$ of the curve singularity $C$ is the
subsemigroup of
$\Z_{{\ge 0}}^r$ which consists of elements of the form
$\vv(g)=(v_1(g), \ldots, v_r(g))$ for all germs $g\in {\cal O}_C$ with
$v_i(g)<\infty$; $i=1, \ldots, r$. The {\bf extended semigroup}
$\widehat
S=\widehat S_C$ of the curve singularity $C$ was defined in \cite{CDG1}.
It is the subsemigroup of $\Z_{{\ge0}}^r\times(\C^*)^r$ which consists
of elements of the form $(\vv(g);\aa(g))=(v_1(g), \ldots, v_r(g);
a_1(g), \ldots, a_r(g))$ for all germs $g\in {\cal O}_C$ with
$v_i(g)<\infty$, $i=1, \ldots, r$ ({\cite{CDG1}}). The extended semigroup
$\widehat S_C$ is well-defined (i.e., does not depend on the choice of
the parameterizations $\varphi_i$) up to a natural equivalence relation.

It is known that both the semigroup $S_C$ and the Alexander polynomial
$\Delta^C(t_1, \ldots, t_r)$ are complete topological invariants of a
plane curve singularity $C$, i.e., each of them determines the germ
$C\subset(\C^2, 0)$ up to topological equivalence ({\cite{W}}, {\cite{Y}}).
The formulae discussed here describe a connection between them.
(In fact from the Eisenbud--Neumann formula for the Alexander polynomial
in terms of a
resolution of a plane curve singularity (see equation~($*$) above) it is
not difficult to understand that the Alexander polynomial $\Delta^C(t_1,
\ldots, t_r)$ may contain with non-zero coefficients only monomials
$\tt^\vv$ for $\vv$ from the semigroup $S_C$ of the curve $C$.)

For a curve singularity $C=\bigcup\limits_{i=1}^{r}C_i\subset (\C^n, 0)$,
let $\pi:\widehat S_C\to \Z^r$ be the natural projection: $(\vv,
\aa)\mapsto \vv$.
For an element $\vv\in \Z^r$, the preimage
$F_{\vv}=\pi^{-1}(\vv)\subset\{\vv\}\times(\C^*)^r\subset\{\vv\}\times\C^r$
is called the fibre of the extended semigroup.
Though $F_{\vv}$ is empty for $\vv\not\in\Z^r_{\ge 0}$, we
define it for
all $\vv\in\Z^r$ in order not to meet (formal) problems with the notations
in the Proof of Theorem~\ref{theo3}.
The fibre $F_{\vv}$ is not empty if and only if $\vv\in S_C$.
For $\vv=(v_1, \ldots, v_r)\in \Z^r$, let $J(\vv)=\{g\in {\cal O}_C:
v_i(g)\ge v_i;\, i=1, \ldots, r\}$ which is an ideal in ${\cal O}_C$. One has a
natural linear map $j_{\vv}:J(\vv)\to\C^r$, which sends $g\in J(\vv)$ to
$(a_1, \ldots, a_r)$, where $a_i$ is the coefficient in the power series
decomposition
$g\circ\varphi_i(\tau_i)=a_i\tau_i^{v_i}+{~terms~of~higher~degree}$
(the number $a_i$ may be equal to zero).
Let $C(\vv)\subset \C^r$ be the image of the map $j_{\vv}$,
$c(\vv)=\mbox{dim }C(\vv)$. It is not difficult to see that
$C(\vv)\cong J(\vv)/J(\vv+{\1})$, where $\1=(1, \ldots, 1)$, and that
$F_{\vv}=C(\vv)\cap(\C^*)^r$ (under the natural
identification of $\{\vv\}\times(\C^*)^r$ and $(\C^*)^r$). Therefore, for
$\vv\in S_C$, the fibre $F_{\vv}$ is the complement to an arrangement of linear
hyperplanes in the linear space $C(\vv)$.

\begin{remark}
For a plane curve singularity $C$,
the extended semigroup $\widehat S_C$ contains some analytic information about
the curve, however the dimensions $c(\vv)$ (and actually the
combinatorial types of the
arrangements of hyperplanes $C(\vv)\cap(\C^r\setminus(\C^*)^r)\subset
C(\vv)$\,) depend only on the topological type of the curve $C$
(see~{\cite{CDG1}}).
\end{remark}

Let ${\cal L}=\Z[[t_1, \ldots, t_r, t_1^{-1}, \ldots, t_r^{-1}]]$ be the
set of formal Laurent series in $t_1, \ldots, t_r$. Elements of ${\cal L}$
are expressions of the form $\sum\limits_{\vv\in\Z^r}k(\vv)\cdot\tt^\vv$ with
$k(\vv)\in\Z$, generally speaking, infinite in all directions. ${\cal L}$
is not a ring, but a $\Z[t_1, \ldots, t_r]$--module (or even a
$\Z[t_1, \ldots, t_r, t_1^{-1}, \ldots, t_r^{-1}]$--module). The polynomial
ring $\Z[t_1, \ldots, t_r]$ can be in a natural way considered as being
embedded into ${\cal L}$.

Let $$L_C(t_1, \ldots, t_r)=\sum\limits_{\vv\in\Z^r}c(\vv)\cdot\tt^\vv\
\in\ {\cal L}.$$
$L_C(\tt)$ is not a power series, but a Laurent series infinitely
long in all directions, since $c(\vv)$ can be positive for $\vv$
with (some) negative components $v_i$ as well. For example, if
there exists a germ $g\in {\cal O}_{\C^n, 0}$ with $v_1(g)=v_1^*$,
then for any $v_2$, \dots, $v_r$ such that $v_i\le v_i(g)$ (including
negative ones), the germ $g$ represents a non-trivial element in
$J(\vv)/J(\vv+\1)$ where $\vv=(v_1^*, v_2, \ldots, v_r)$,
$J(\vv)=\{g\in{\cal O}_{\C^n, 0}: \vv(g)\ge\vv\}$. One can understand
that along each line in the lattice $\Z^r$ parallel to a coordinate
one the coefficients $c(\vv)$ stabilize in each direction, i.e., if
$v_i'$ and $v_i''$ are negative, or if $v_i'$ and $v_i''$ are positive
and large enough, then
$c(v_1,\ldots, v_i',\ldots, v_r) = c(v_1,\ldots, v_i'',\ldots, v_r)$.
This implies that
$$P'_C(t_1, \cdots ,t_r)= L_C(t_1, \cdots, t_r)\cdot
\prod\limits_{i=1}^r (t_i-1)$$
is a polynomial (it also follows from the proof of Theorem~\ref{theo3}).

\begin{proposition}\label{prop1}
For a curve singularity $C=\bigcup\limits_{i=1}^{r}C_i\subset (\C^n, 0)$,
$r>1$, the polynomial $P^\prime_C(t_1, \ldots, t_r)$ is divisible by
$(t_1\cdot\ldots\cdot t_r-1)$, i.e., the power series
$$
P_C(t_1, \ldots, t_r)=P^\prime_C(t_1, \ldots,
t_r)/(t_1\cdot\ldots\cdot t_r-1) \in \Z[[t_1, \ldots, t_r]]
$$
is, in fact, a polynomial.
\end{proposition}

The proposition follows from the proof of Theorem~\ref{theo3}
(see below).

\smallskip
We call $P_C(t_1, \ldots, t_r)$, $r>1$, the (generalized)
{\bf Poincar\'e
polynomial} of the curve singularity $C$. For $r=1$, $P_C(t)$ is
not
a polynomial, but a power series and it coincides with usual Poincar\'e
series of the filtration in $\cal O_C$ defined by a normalization.

\begin{remark}
Theorem~\ref{theo1} implies that, for a plane curve singularity,
the polynomial $P_C(\tt)$ determines the semigroup $S_C$. This is
not the case for non--plane curves (see~{\cite{CDK}}).
\end{remark}

\subsection{The semigroup of an irreducible plane curve
singularity.}\label{irred}
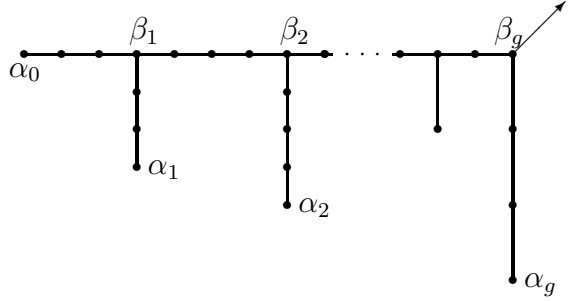
\begin{figure}
$$
\unitlength=1.00mm
\begin{picture}(80.00,30.00)(-10,3)
\thicklines
\put(-5,30){\line(1,0){41}}
\put(44,30){\line(1,0){16}}
\put(38,30){\circle*{0.5}}
\put(40,30){\circle*{0.5}}
\put(42,30){\circle*{0.5}}
\put(30,10){\line(0,1){20}}
\put(50,20){\line(0,1){10}}
\put(60,0){\line(0,1){30}}
\put(10,15){\line(0,1){15}}
\thinlines
\put(60,30){\vector(1,1){7}}
\put(20,30){\circle*{1}}
\put(30,30){\circle*{1}}
\put(50,30){\circle*{1}}
\put(60,30){\circle*{1}}
\put(30,20){\circle*{1}}
\put(60,20){\circle*{1}}
\put(60,10){\circle*{1}}
\put(10,30){\circle*{1}}
\put(30,10){\circle*{1}}
\put(50,20){\circle*{1}}
\put(60,0){\circle*{1}}
\put(-5,30){\circle*{1}}
\put(0,30){\circle*{1}}
\put(5,30){\circle*{1}}
\put(15,30){\circle*{1}}
\put(25,30){\circle*{1}}
\put(35,30){\circle*{1}}
\put(45,30){\circle*{1}}
\put(55,30){\circle*{1}}
\put(10,25){\circle*{1}}
\put(10,20){\circle*{1}}
\put(10,15){\circle*{1}}
\put(30,25){\circle*{1}}
\put(30,15){\circle*{1}}
\put(35,30){\circle*{1}}
\put(-7,27){$\alpha_0$}
\put(11.5,14){$\alpha_1$}
\put(31.5,9){$\alpha_2$}
\put(61.5,-1){$\alpha_g$}
\put(9,32){$\beta_1$}
\put(29,32){$\beta_2$}
\put(57.5,32){$\beta_g$}
\end{picture}
$$
\caption{The dual resolution graph of the curve $C$.}
\label{fig0}
\end{figure}

If the curve $C$ is irreducible its dual graph looks like
on Fig.~\ref{fig0}.
After Zariski (see e.g. \cite{ZT}) it is known that the set of
elements $\{\betab_j :=m^{\alpha_j} : 0 \le  j\le g\}$ is the
minimal system of generators of the semigroup of values
$S_{C}\subset \Z_{\ge 0}$; moreover $m^{\beta_j} =
(n_j+1)m^{\alpha_j}$ for some integers $n_j$, $j=1,\ldots,g$
(integers $n_j+1$ are in fact parts of the Puiseux pairs of the
curve; the vertices $\alpha_j$ and $\beta_j$ are indicated on
Fig.~\ref{fig0}). In what follows we shall use the following properties of
the minimal embedded resolution and of the semigroup $S_C$ of
an irreducible curve singularity (see, e.g., \cite{ZT}, \cite{D}):
\begin{enumerate}
\item There is only a finite number of positive integers which
do not belong to the semigroup $S_C$ and the largest one $\delta-1$
is equal to $\sum_{j=1}^g n_j \betab_j- \betab_0$ ($\delta$ is called the
conductor of the semigroup $S_C$).
\item Each element $v\in S_{C}$ can be uniquely represented in the form
$v = k_0 \betab_0 + \sum_{j=1}^{g}k_j\betab_j$ with
$k_0\ge 0$ and $0\le k_j\le n_j$ for $1 \le j \le g$.
\item $(n_j+1)\betab_{j} <\betab_{j+1}$ for $j=1,\ldots, g-1$.
\item The element $(n_j+1)\betab_{j}$ belongs to the semigroup
$\langle\betab_0,\ldots,\betab_{j-1}\rangle$ generated by
$\betab_0$, \dots, $\betab_{j-1}$ ($j=1,\ldots, g$).
\item If, for a germ $\varphi\in {\cal O}_{\C^2,0}$, the strict
transform of the curve $\{\varphi=0\}$ intersects only divisors
$\sigma$ with $\sigma < \beta_j$
(i.e., those which lie between $\alpha_0$ and $\beta_j$) then
$v(\varphi)\in \langle \betab_0, \ldots ,\betab_{j-1}\rangle$.
\end{enumerate}

According to {\cite{A'C}} the zeta-function $\zeta_C(t)$ of the
singularity $C=\{f=0\}$ is equal to
$\prod\limits_\sigma(1-t^{v_\sigma})^{-\chi({\stackrel{\circ}{E}}_\sigma)}$.
The property 2 (unique decomposition) permits to easily prove the
equality of the Poincar\'e series $P_C(t)$ and the zeta-function
$\zeta_C(t)$ in this case (\cite{CDG2}; in fact the proof is almost
repeated in the proof of Proposition~\ref{prop3},
subsection~\ref{s-reduction}).

\subsection{Graded topological spaces and the Euler
characteristic.}\label{s-graded}
To describe some constructions in the proof of Theorems~\ref{theo1} and
\ref{theo2}, it is convenient to use the notion of a graded space.
By a graded space (with $r$--grading) we shall have in mind
a disjoint union $Z$ of
topological spaces $Z_\vv$ corresponding to elements $\vv$
from $\Z_{\ge 0}^r$.
We shall write $Z=\sum\limits_{\vv\in\Z_{\ge 0}^r}Z_\vv\cdot\tt^\vv$.
If $Z'$ and $Z''$ are graded spaces, then their sum (disjoint union)
$Z' + Z''$ and their product $Z'\times Z''$ are graded spaces as well
(in the natural sense;
e.g., if $Z'=\sum\limits_{\vv\in\Z_{\ge 0}^r}Z_\vv '\cdot\tt^\vv$,
$Z''=\sum\limits_{\vv\in\Z_{\ge 0}^r}Z_\vv ''\cdot\tt^\vv$, then
$Z'\times Z'' = \sum\limits_{\vv\in\Z_{\ge 0}^r}
\sum\limits_{\vv'+\vv''=\vv}(Z_{\vv'}'\times Z_{\vv''}'')\cdot\tt^\vv$).
A map $Z'\to Z''$ of graded spaces is a set of maps
$Z'_\vv\to Z''_\vv$. For a graded space $Z=\sum\limits_{\vv\in\Z_{\ge 0}^r}
Z_\vv\cdot\tt^\vv$ its Euler characteristic (see the remarks below)
$\chi(Z)$ is the series
$\sum\limits_{\vv\in\Z_{\ge 0}^r}\chi(Z_\vv)\cdot\tt^\vv\in
\Z[[t_1, \ldots, t_r]]$. One has $\chi(Z' + Z'')=\chi(Z')+\chi(Z'')$,
$\chi(Z'\times
Z'')=\chi(Z')\cdot\chi(Z'')$. A graded semigroup is a graded
space with a (commutative) semigroup operation which respects the
grading. The extended semigroup of a plane curve singularity is
in the natural sense a graded semigroup.

There exist somewhat different definitions of the Euler characteristic
which do not coincide for non-compact sets. (For compact spaces
(say, for projective varieties) all definitions of the Euler
characteristic are essentially the same.)
The most usual definition of the Euler characteristic of a
topological space (say, of a $CW$--complex) $X$ is $\chi(X)=
\sum\limits_{q\ge 0}(-1)^q\dim H_q(X;\R)$. If $X=X_1\cup X_2$ where the
spaces ($CW$--complexes) $X$, $X_1$ and $X_2$ are compact, one has
$\chi(X)=\chi(X_1)+\chi(X_2)-\chi(X_1\cap X_2)$. Therefore for compact
spaces the Euler characteristic possesses the additivity property.
This permits to consider it as a generalized (nonpositive) measure on
the algebra of such spaces. However spaces we are interested in (e.g.,
fibres of the extended semigroup of a curve) are noncompact semialgebraic
sets (complex or real). The Euler characteristic defined above does not
possess the additivity property for such spaces. For example, let $X$ be the
circle $S^1$, let $X_1$ be a point of $X$, and let $X_2=X\setminus X_1$ be
a (real) line. Then one has $\chi(X)=0$, $\chi(X_1)=\chi(X_2)=1$,
$\chi(X_1\cap X_2)=\chi(\emptyset)=0$, and $0\ne 1+1-0$.

In order to have the desired additivity property one should define the
Euler characteristic $\chi(X)$ of a semialgebraic space $X$
(the difference of two projective spaces) as
$$\sum\limits_{q\ge 0}(-1)^q\dim H_q(X^*, *;\R),$$ where $X^*$ is the
one-point compactification of the space $X$ (if $X$ is compact, the
one-point compactification of it is the disjoint union of $X$ with a
point), $*$ is the added ("infinite") point. We shall use this definition.
The algebra generated by semialgebraic sets consists of
constructible
sets. A constructible set can be represented as the disjoint union
of a finite number of semianalytic sets. The Euler characteristic of
a constructible set should be defined as the sum of Euler
characteristics of the corresponding semialgebraic sets.
One can show that the Euler characteristic defined this way does possess
the additivity property (in the example above
$\chi(X)=0$, $\chi(X_1)=1$, $\chi(X_2)=-1$). Moreover, a constructible
set $X$ can be represented as a disjoint union of {\bf a finite number}
of open cells so that the boundary of a cell of some dimension (its
closure in $X$ minus itself) lies in the union of cells of smaller dimensions.
(This does not mean a representation of the space $X$ as a
$CW$--complex since in general (for noncompact sets) one does not have
maps of closed balls into $X$ which determine the cells. For example the
real line $\R^1$ is simply one cell of dimension 1.) One can see that the
Euler characteristic of the (constructible) set $X$ is equal to the
alternating sum of numbers of cells of different dimensions. The Euler
characteristic also possesses the multiplicativity property:
$\chi(X_1\times X_2)=\chi(X_1)\cdot\chi(X_2)$.

\medskip
Let $X$ be a topological space and
let $S^kX$ be its $k$-th symmetric
power $($$S^kX=\underbrace{X\times\ldots\times X}_{\mbox{$k$ times}}/S_k$,
i.e., the space of unordered $k$--tuples of points of $X$, where $S_k$ is the
group of permutations of $k$ elements; $S^0X=\bullet$ is a point$)$.
The graded space (with $1$-granding) $S^0X + S^1X \cdot t+
S^2X\cdot t^2 + \cdots$ is in
fact a graded semigroup with the semigroup operation defined by the
union of $k$--tuples of points.

\begin{lemma}\label{lemma1}
$$
\chi (S^0X + S^1X \cdot t+ S^2X\cdot t^2 +\cdots) =
(1-t)^{-\chi(X)}.
$$
\end{lemma}

\begin{proof}
Let us denote $\chi(X)$ simply by $\chi$. The coefficient at $t^k$
in $(1-t)^{-\chi}$ is equal to
$$
{\chi+k-1\choose k} = (-1)^k {-\chi\choose k} =
\frac{\chi (\chi+1)\cdot\ldots\cdot (\chi+k-1)}{1\cdot\ldots\cdot k}.
$$

It is clear that $\chi(S^kX)$ is a polynomial in $\chi$.
To show that this polynomial coincides with ${\chi+k-1\choose k}$
one can check it for an infinite set of values of $\chi$.
Let us take $\chi$ positive, and let $X$ be the (disjoint) union of
$\chi$ points ($\chi(X)=\chi$). In this case $S^k X$ consists of
${\chi+k-1\choose k}$ points and thus $\chi (S^kX) = \# S^k X
= {\chi+k-1\choose k}$.
\end{proof}

\begin{remark}
For $X =\{Z_1, \ldots, Z_\chi\}$ ($Z_i$ are points), this formula
is a consequence of the equation
$$
\left(1\,+\,\sum_{k=1}^{\infty} (\sum_{n_1+\cdots +n_\chi=k}
\underline{Z}^{\underline{n}})t^k\right)\cdot\prod_{i=1}^{\chi}(1-Z_i t)
\,=\,1.
$$
\end{remark}

\begin{corollary}
If $\chi(X)\le 0$ and $k\ge -\chi(X)+1$, then $\chi(S^kX)=0$.
In particular, for $X=\CP^1\setminus\{\mbox{$s$ points}\}$,
$s\ge 2$, $\chi(S^kX)=0$ for $k\ge s-1$.
\end{corollary}

The group $\C^*$ of non-zero complex numbers acts (freely) on
$\Z_{{\ge0}}^r\times(\C^*)^r$ (by multiplication of all the
coordinates $a_i$). The corresponding factor--space
$\Z_{{\ge0}}^r\times(\C^*)^r/\C^*=\Z_{{\ge0}}^r\times\P((\C^*)^r)=
\sum\limits_{\vv\in\Z_{\ge 0}^r}\P((\C^*)^r)\cdot\tt^\vv$ has a
natural structure of a semigroup. The extended semigroup
$\widehat S_C\subset\Z_{{\ge0}}^r\times(\C^*)^r$ is invariant with respect to
the $\C^*$--action. The factor--space $\P\widehat S_C=\widehat S_C/\C^*$
will be called the projectivization of the extended semigroup of the
curve $C$ (it is also a graded semigroup in the natural sense). One has
$\P\widehat S_C=\sum\limits_{\vv\in\Z_{\ge 0}^r}\P F_{\vv}\cdot\tt^\vv$,
where $\P F_{\vv}=F_{\vv}/\C^*$
is the projectivization of the fibre $F_{\vv}$. For $\vv\in S_C$, the space
$\P F_{\vv}$ is the complement to an arrangement of projective
hyperplanes in a $(c(\vv)-1)$--dimensional complex projective space
$\P C(\vv)$, which can be identified with the complement to an arrangement
of affine hyperplanes in a complex affine space after choosing one of the
coordinates $a_i$. Let ${\underline{\delta}}=(\delta_1, \ldots, \delta_r)$
be the {\bf conductor} of the semigroup $S_C$ of the curve $C$, i.e.,
the minimal element for which $\underline{\delta} + \Z^r_{\ge 0}\subset S_C$.
If $\vv\ge {\underline{\delta}}$ (i.e., if $v_i\ge\delta_i$
for all $i=1, \ldots, r$), then the fibre $F_{\vv}$ of the
extended semigroup coincides with $(\C^*)^r$ and the Euler characteristic
$\chi(\P F_{\vv})$ of its projectivization is equal to $0$
(for $r>1$; for $r=1$ it is equal to $1$). Moreover one can show
that the Euler characteristic
$\chi(\P\widehat S_C)$ is a polynomial in $t_1$, \dots, $t_r$ (see
Theorem~\ref{theo3}). This polynomial participates in the formulation
of Theorem~\ref{theo2}.

\section{Proofs}

The following statement shows that Theorems 1 and 2 are
equivalent to each other.

\begin{theorem}\label{theo3}
For an arbitrary (i.e., not necessarily plane) curve singularity
$C=\bigcup\limits_{i=1}^{r}C_i\subset (\C^n, 0)$, $r>1$,
$\chi(\P\widehat S_C)$ is a polynomial and $\chi({\Bbb P}\widehat
S_C)\cdot(t_1\cdot\ldots\cdot t_r-1)=P^\prime_C(\tt)$. As a
consequence, $P_C(\tt)$ is a polynomial and
$\chi({\Bbb P}\widehat S_C)=P_C(\tt)$.
\end{theorem}

\begin{proof}
Let $\ww$ be an element of $\Z_{\ge 0}^r$, and let $b_\vv=
\dim{J(\vv)/J(\ww)}$. For $I\subset I_0=\{1, 2, \ldots, r\}$,
let $\#I$ be the number of elements in $I$, and let $\1_I$ be
the element of $\Z_{\ge 0}^r$, the $i$-th component of which is equal
to $1$ (respectively to $0$) if $i\in I$ (respectively if $i\not\in I$).
One has $\1_{I_0}=\1$. Let $L_I\subset\C^r$ be the subspace
$\{(a_1, \ldots, a_r)\in\C^r: a_i=0 \mbox{ for } i\in I\}$.

One has
$$
\chi(\P F_{\vv})=\chi(\P C(\vv)) - \chi(\bigcup\limits_{i=1}^r
\P(C(\vv)\cap L_{\{i\}}))=
$$
$$
=\chi(\P C(\vv)) - \sum\limits_{I\subset I_0, I\ne\emptyset}
(-1)^{\# I-1}\chi(\P(C(\vv)\cap L_{I}))=
$$
$$
=
\sum\limits_{I\subset I_0}(-1)^{\# I}\chi(\P(C(\vv)\cap L_{I}))=
\sum\limits_{I\subset I_0}
(-1)^{\# I}\dim(C(\vv)\cap L_{I}).
$$
If $\vv\le\ww-\1$, $\dim(C(\vv)\cap L_{I})=b_{\vv+\1_I}-b_{\vv+\1}$
and therefore $\chi(\P F_{\vv})=\sum\limits_{I\subset I_0}(-1)^{\#I}
(b_{\vv+\1_I}-b_{\vv+\1})$.
This implies that the coefficient at $\tt^\vv$
in the series $\chi({\Bbb P}\widehat S_C)\cdot(t_1\cdot\ldots\cdot t_r-1)$
is equal to
$$
\sum\limits_{I\subset I_0}(-1)^{\#I}(b_{\vv+\1_I-\1}-b_{\vv})-
\sum\limits_{I\subset
I_0}(-1)^{\#I}(b_{\vv+\1_I}-b_{\vv+\1})
$$
and, since $\sum_{I\subset I_{0}}(-1)^{\#I}=0$, also to
$$
\sum\limits_{I\subset
I_0}(-1)^{\#I}(b_{\vv-\1+\1_I}-b_{\vv+\1_I})=
\sum\limits_{I\subset I_0}(-1)^{\#I}c(\vv-\1+\1_I).
$$
The coefficient at $\tt^\vv$ in the polynomial $P^\prime_C(\tt)=
(\sum c(\vv)\tt^\vv)\cdot\left(\prod\limits_{i=1}^r(t_i-1)\right)$
is also equal to $\sum\limits_{I\subset I_0}(-1)^{\#I}c(\vv-\1+\1_I)$.

The facts that the series $\chi({\Bbb P}\widehat S_C)$ does not
contain (with non-zero coefficients) monomials $\tt^\vv$ with
$\vv\ge {\underline{\delta}}$ and $P^{\prime}_C(\tt)$ is a
polynomial imply that $\chi({\Bbb P}\widehat S_C)$ is a polynomial as
well.
\end{proof}

\smallskip\noindent{\bf Proof of Theorem 2.}
By the Eisenbud-Neumann formula, Theorem~\ref{theo2} follows from

\begin{theorem}\label{theo4}
$\chi(\P\widehat S_C)=\prod\limits_{\sigma}
(1-\tt^{\mm^\sigma})^{-\chi({\stackrel{\circ}{E}}_\sigma)}$.
\end{theorem}

The main course of {\bf the proof} for Theorem 4 goes as follows.
We use an embedded resolution of the curve $C$. In terms of it
in \ref{s-spaceY} we construct a graded space (with $r$--grading;
in fact a graded semigroup) $Y$ such that its Euler characteristic
is equal to the Alexander polynomial $\Delta^C(t_1, \ldots, t_r)$.
The construction of the space $Y$ is
natural after Lemma \ref{lemma1} above which provides the
equality of the Euler characteristic of $Y$ and the Alexander
polynomial $\Delta^C(\tt)$. Points of $Y$ are
represented by unordered $k$--tuples (with different $k$) of
points of the exceptional divisor of the resolution without the
self-intersection points of it and the intersection points of it
with the strict transform of the curve $C$. The semigroup
operation is defined by the union of $k$--tuples. The grading is
defined by the multiplicities $\mm^\sigma$ of the components
$E_\sigma$ of the exceptional divisor on which the points lie.

We construct
a map (a graded semigroup homomorphism) $\Pi$ from $Y$ to the
projectivisation $\P \widehat S_C$ of the extended semigroup $\widehat S_C$ which is
surjective up to a grading high enough, i.e., $\Pi$ maps $Y_\vv$ onto
the fibre $\P F_{\vv}$ for $\vv \le V$, where $V$ is an arbitrary
point of $\Z_{\ge 0}^r$ chosen in advance. This map is defined
in the following way. For a point $y\in Y$, we take a function
$g\in{\cal O}_{\C^2, 0}$ such that the strict transform of the curve
$\{g=0\}$ intersects the exceptional divisor of the resolution of the
curve $C$ just at the points which define the element $y$ and we
put $\Pi(y)=(v_1(g), \ldots, v_r(g); a_1(g):\ldots:a_r(g))$.

If $\Pi$ would be injective this would complete the proof.
However the space $Y$ is too big and $\Pi$ is very far from being
injective. In \ref{s-reduction} we reduce $Y$ to another space,
$\widetilde{Y}$, together with the corresponding map
$\widetilde{\Pi}: \widetilde{Y} \to \P (\widehat{S}_C)$ so that
$\chi(\widetilde{Y})=\chi(Y)=\Delta^C(\tt)$, $\widetilde{\Pi}(\widetilde Y)=
\Pi(Y)$. Roughly speaking, this reduction consists of excluding all
dead ends of the resolution graph. This is a way to omit some
obvious repeated images by $\widetilde{\Pi}$.
This permits to pay the main attention only to those components
$E_\sigma$, $\sigma\in\Gamma$, of the exceptional divisor $D$,
for which the number $s_\sigma$ of essential points is $\ge 2$.
In Lemma \ref{lemma3} we show that, if
$\widetilde{\Pi}(y_1)=\widetilde{\Pi}(y_2)$
($y_i \in \widetilde{Y}$), then $y_1$ and $y_2$ have "many"
points ($\ge s_\sigma-1\ge 1$) on some components $E_\sigma$
of the exceptional divisor (generally speaking, different for
$y_1$ and $y_2$). If an element $y\in\widetilde{Y}$
is represented by a $k$-tuple with at least $s_\sigma-1$ of
them on the component $E_\sigma$ of the exceptional divisor
with $s_\sigma \ge 2$, we say that $\sigma$ is a {\bf cut}
of $y$ (or rather of the connected component of $\widetilde{Y}$
where $y$ lies).

Moreover, in \ref{s-cuts} we prove more fine statements about
the distribution of the cuts of $y_1$ and $y_2$ in this case.
Namely, up to the numbering of $y_1$ and $y_2$ there exists a cut
$\sigma$ of $y_1$ such that for each strict transform
$\widetilde C_i$ of a branch $C_i$ of the curve $C$ greater than $\sigma$
there is a cut of $y_2$ on the geodesic from $\sigma$
to $\widetilde C_i$ on the dual graph of the resolution.
This is the most complicated (combinatorial) part of the proof.

To prove that $\chi(\mbox{\it Im\,}\widetilde{\Pi})=\chi(\widetilde Y)$,
in \ref{s-torus} we analyse places where the map $\widetilde{\Pi}$ is not
injective. At such places we indicate some parts of the space
$\widetilde Y$ which are fibred into complex tori $(\C^*)^{s-1}$,
$s\ge 2$ (and thus have zero Euler characteristic). Removing these
parts does not change the Euler characteristic of the source
and does not change the image $\widetilde{\Pi}(\widetilde Y)$.
This permits to obtain a one-to-one correspondence
between them. This completes the proof.

One can say that the proof consists of an explicit computation
of $\chi(\P \widehat S_C)$ and of $\Delta^C(t_1, \ldots, t^r)$
in terms of an embedded resolution of the curve $C$ which shows
that they coincide. At the moment a direct proof which explains
why the Euler characteristic $\chi(\P \widehat S_C)$ coincides
with the Alexander polynomial is not known.

\subsection{Construction of the space $Y$ and the map $\Pi$.}
\label{s-spaceY}
One can easily see that the right--hand side of the equation of
Theorem~\ref{theo4} does not depend on the resolution of the curve $C$.
Let $\underline{V}=(V_1,\ldots, V_r)$ be an arbitrary point of the
lattice $\Z_{\ge0}^r$. Let us take the minimal (embedded) resolution
of the curve $C$ and let us make additional blow--ups {\bf of intersection
points} of components of the total transform of the curve $C$ so that,
for each function $g\in{\cal O}_{\C^2, 0}$ with $\vv(g)\le\underline{V}$,
the strict transform of the curve $\{g=0\}$ intersects the total
transform $(f\circ\pi)^{-1}(0)$ of the curve $C$ only at smooth
points of $(f\circ\pi)^{-1}(0)$. (These additional blow-ups do
not change dead ends and star points of the dual graph $\Gamma$
of the resolution. Moreover, one can see that, if $\underline V$
is big enough (say, $\underline{V}\ge\underline{\delta}$), then
each strict transform $\widetilde C_i$ lies on its individual
"long" branch of the graph $\Gamma$.) We fix such a
 resolution for the rest of the paper. Let
$$
Y=\prod\limits_{\sigma}(\bullet+
S^1{\stackrel{\circ}{E}}_\sigma\cdot\tt^{\mm^\sigma}
+S^2{\stackrel{\circ}{E}}_\sigma\cdot\tt^{2\mm^\sigma}+\ldots)
$$
where
$S^k{\stackrel{\circ}{E}}_\sigma$ is the $k$th symmetric power of
${\stackrel{\circ}{E}}_\sigma$ and $\bullet$ is a point
($=S^0{\stackrel{\circ}{E}}_\sigma$).
Lemma~\ref{lemma1} implies that
$\chi(Y)=\prod\limits_{\sigma}
(1-\tt^{\mm^\sigma})^{-\chi({\stackrel{\circ}{E}}_\sigma)}$.

Let us define a map $\Pi:Y\to\P\widehat S_C$ as follows. One has
$$
Y=\sum\limits_{\{k_\sigma\}}\left(\prod\limits_\sigma
S^{k_\sigma}{\stackrel{\circ}{E}}_\sigma
\cdot\tt^{\sum k_\sigma\mm^\sigma}\right).
$$
A point $y$ of the space
$\prod\limits_\sigma (S^{k_\sigma}{\stackrel{\circ}{E}}_\sigma\cdot
\tt^{k_\sigma\mm^\sigma})$ is represented by a set of smooth
points of the exceptional divisor $D$ (i.e., of
${\stackrel{\circ}{D}}=
\bigcup\limits_{\sigma}{\stackrel{\circ}{E}}_\sigma$) with $k_\sigma$ points
$Q_1^\sigma$, \dots, $Q_{k_\sigma}^\sigma$ on the component
${\stackrel{\circ}{E}}_\sigma$. For a point $A\in{\stackrel{\circ}{D}}$,
let ${\widetilde L}_A$ be a germ of a nonsingular (complex
analytic) curve
transversal to the exceptional divisor $D$ at the point $A$. Let the image
$L_A=\pi({\widetilde L}_A)\subset (\C^2, 0)$ of the curve ${\widetilde L}_A$
be given by an equation $\{g_A=0\}$ ($g_A\in{\cal O}_{\C^2, 0}$). By
definition $\Pi(y)\in\P\widehat S_C$ is represented by the element
$(\vv(g), \aa(g))\in\widehat S_C$, where
$g=\prod\limits_{\sigma}\prod\limits_{j=1}^{k_\sigma}g_{Q_j^\sigma}$.

\begin{lemma}\label{lemma2}
The element $\Pi(y)\in\P\widehat S_C$ does not depend on the
choice of curves ${\widetilde L}_A$.
\end{lemma}

\begin{proof}
Let ${\widetilde L}_A^\prime$ be another germ of a nonsingular
(complex) curve transversal to the exceptional divisor $D$ at
the point $A\in{\stackrel{\circ}{D}}$,
$L_A^\prime=\pi({\widetilde L}_A^\prime)=\{g_A^\prime=0\}$, and let
$g^\prime=\prod\limits_{\sigma}\prod\limits_{j=1}^{k_\sigma}g_{Q_j^\sigma}^\prime$.
Let ${\widetilde g}=g\circ\pi$ and ${\widetilde g}^\prime=g^\prime\circ\pi$
be the
liftings of the functions $g$ and $g^\prime$ to the space $X$ of the
resolution,
and let $\psi=\widetilde g^\prime/\widetilde g$ be their ratio.
The function $\psi$ has zeros along
the curves ${\widetilde L}_{Q_j^\sigma}^\prime$ and poles along the curves
${\widetilde L}_{Q_j^\sigma}$. Therefore the restriction of the function
$\psi$ to the exceptional divisor $D$ is a regular (holomorphic) function
on $D$ and thus $\psi$ is a constant (say, $c$) on $D$. It implies that
$\vv(g^\prime)=\vv(g)$, $\aa(g^\prime)=c\cdot\aa(g)$ and
therefore the elements in
$\P\widehat S_C$, represented by $(\vv(g), \aa(g))$ and $(\vv(g^\prime),
\aa(g^\prime))$,
coincide.
\end{proof}

\begin{remark}
In fact one can say that $Y$ is a graded semigroup (with respect
to the operation defined by the union of sets) and $\Pi$ is a graded
semigroup homomorphism. Then it is sufficient to define $\Pi$ only
for monomials of the form $[A]\cdot \tt^{\mm^\sigma}$, where $A$ is
a point of ${\stackrel{\circ}{E}}_\sigma$. We shall use this to
define the map $\widetilde\Pi$ below.
\end{remark}

\begin{proposition}\label{prop2}
For $\vv\le\underline{V}$ one has $(Im\,\Pi)_\vv=\P F_\vv$.
\end{proposition}

\begin{proof}
If $g\in{\cal O}_{\C^2, 0}$ and $\vv(g)\le{\underline{V}}$, then
the strict transform $\widetilde L$ of the curve $\{g=0\}$
intersects the exceptional divisor $D$ only at
smooth points of $(f\circ\pi)^{-1}(0)$. Let $Q^\sigma_j$ ($j=1, \ldots,
k_\sigma$) be the points of
intersection of the curve $\widetilde L$ with the component
$E_\sigma$ of the exceptional divisor counting with their
multiplicities, i.e., each point is taken as many times as the
intersection number of the curve $\widetilde L$ with the
component $E_\sigma$ at it, let $K_\sigma=\{Q_1^\sigma, \ldots,
Q_{k_\sigma}^\sigma\}\subset{\stackrel{\circ}{E}}_\sigma$.
For a subset $K_\sigma$ of ${\stackrel{\circ}{E}}_\sigma$
(or of $\widetilde E_\sigma$) with $\# K_\sigma=k_\sigma$,
by $[K_\sigma]$ we shall denote the corresponding point of the
space $S^{k_\sigma}{\stackrel{\circ}{E}}_\sigma$ (or of
$S^{k_\sigma}\widetilde E_\sigma$). Then one can see that
$$
(\vv(g),
\aa(g))=\left(
\vv(\,\prod\limits_\sigma\prod\limits_{j=1}^{k_\sigma}g_{Q_j^\sigma}),\,
\aa(\,\prod\limits_\sigma\prod\limits_{j=1}^{k_\sigma}g_{Q_j^\sigma})
\right)=
\Pi\left(\prod\limits_\sigma
([K_\sigma]\cdot\tt^{k_\sigma\mm^\sigma})\right)
$$ (the proof repeats the one of Lemma~\ref{lemma2}).
\end{proof}

\subsection{Reduction of the graded space
$Y$.}\label{s-reduction}
In order to reduce the space $Y$ we will use arithmetical
properties of the semigroup $S_{C}$ and its relation with the
dual graph of a resolution which can be found in \cite{D1}, section (3.20)
(see also \cite{D}, section 1).

Let $D^\prime$ be the union of components $E_\sigma$ with
at least two essential points, i.e., with $s_\sigma\ge 2$, and let
$\Delta^\prime$ be the set of the corresponding vertices.
Connected components of the complement
$D\setminus D^\prime$, which do not contain the starting divisor
${\bf 1}$, are tails of the dual graph $\Gamma$ of the
resolution and correspond to (some) dead ends $\delta$ of the
graph $\Gamma$. Let $\Delta$ be the set of these dead ends. For
$\delta\in\Delta$, $st_\delta$ is the vertex of $\Delta^\prime$
such that $E_{st_\delta}$ intersects the corresponding connected
component of $\overline{D\setminus D^\prime}$. Let
$$
Y^\prime=\prod\limits_{\sigma\in\Delta^\prime}(\bullet+
S^1{{\widetilde{E}}_\sigma}\cdot\tt^{\mm^\sigma}
+S^2{{\widetilde{E}}_\sigma}\cdot\tt^{2\mm^\sigma}+\ldots).
$$
Pay attention that the spaces ${\stackrel{\circ}{E}}_\sigma$
in the definition of $Y$ are substituted here by the spaces
${\widetilde E}_\sigma$. The possibility to deal with points
of $\widetilde E_\sigma \setminus{\stackrel{\circ}{E}}_\sigma$
in the same way as with other points of $\widetilde E_\sigma$
will be explained below (in the Remark after Proposition~\ref{prop4}).

For a dead end $\delta\in\Delta$, $\mm^{st_\delta}$ is a multiple of
$\mm^{\delta}$: $\mm^{st_\delta}=(n_\delta+1)\cdot \mm^{\delta}$
(the number $n_{\delta}$ is the corresponding $n_j$ which appears
in  \ref{irred} for a branch $C_i$ such that the dead end
$\delta$
belongs to the minimal resolution graph of the curve $C_i$). Let
$Y_\delta=\sum\limits_{k=0}^{n_\delta}\bullet\cdot\tt^{k\mm^\delta}$.

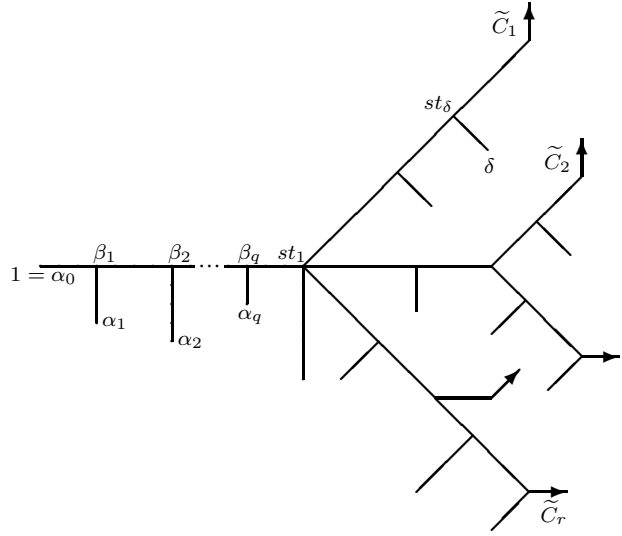
\begin{figure}
$$
\unitlength=0.50mm
\begin{picture}(120.00,100.00)(10,-40)
\thicklines
\put(-10,30){\line(1,0){41}}
\put(39,30){\line(1,0){21}}
\put(60,30){\line(1,0){50}}
\put(90,30){\line(0,-1){12}}

\put(60,30){\line(1,1){60}}
\put(120,90){\vector(0,1){10}}
\put(110,92){{\scriptsize $\widetilde C_1$}}
\put(85,55){\line(1,-1){9}}
\put(100,70){\line(1,-1){9}}
\put(108,55){{\scriptsize $\delta$}}
\put(92,72){{\scriptsize $st_\delta$}}

\put(60,30){\line(1,-1){60}}
\put(120,-30){\vector(1,0){10}}
\put(123,-38){{\scriptsize $\widetilde C_r$}}
\put(80,10){\line(-1,-1){10}}
\put(105,-15){\line(-1,-1){15}}
\put(120,-30){\line(-1,-1){10}}
\put(95,-5){\line(1,0){15}}
\put(110,-5){\vector(1,1){7.5}}

\put(110,30){\line(1,1){24}}
\put(122,42){\line(1,-1){9}}
\put(134,54){\vector(0,1){10}}
\put(124,56){{\scriptsize $\widetilde C_2$}}

\put(110,30){\line(1,-1){24}}
\put(119,21){\line(-1,-1){9}}
\put(134,6){\line(-1,-1){9}}
\put(134,6){\vector(1,0){10}}

\put(33,30){\circle*{0.5}}
\put(35,30){\circle*{0.5}}
\put(37,30){\circle*{0.5}}
\put(25,10){\line(0,1){20}}
\put(45,20){\line(0,1){10}}
\put(60,0){\line(0,1){30}}
\put(5,15){\line(0,1){15}}
\thinlines
\put(15,30){\circle*{1}}
\put(25,30){\circle*{1}}
\put(45,30){\circle*{1}}
\put(60,30){\circle*{1}}
\put(25,20){\circle*{1}}
\put(60,20){\circle*{1}}
\put(60,10){\circle*{1}}
\put(5,30){\circle*{1}}
\put(25,10){\circle*{1}}
\put(45,20){\circle*{1}}
\put(60,0){\circle*{1}}

\put(-10,30){\circle*{1}}
\put(-5,30){\circle*{1}}
\put(0,30){\circle*{1}}
\put(10,30){\circle*{1}}
\put(20,30){\circle*{1}}
\put(30,30){\circle*{1}}
\put(40,30){\circle*{1}}
\put(50,30){\circle*{1}}
\put(55,30){\circle*{1}}
\put(60,25){\circle*{1}}
\put(60,15){\circle*{1}}
\put(60,5){\circle*{1}}
\put(5,25){\circle*{1}}
\put(5,20){\circle*{1}}
\put(5,15){\circle*{1}}
\put(25,25){\circle*{1}}
\put(25,15){\circle*{1}}
\put(-18,26){{\scriptsize $1=\alpha_0$}}
\put(6.5,14){{\scriptsize$\alpha_1$}}
\put(26.5,9){{\scriptsize$\alpha_2$}}
\put(42.5,16){{\scriptsize$\alpha_q$}}
\put(4,32){{\scriptsize$\beta_1$}}
\put(24,32){{\scriptsize$\beta_2$}}
\put(42.5,32){\scriptsize{$\beta_q$}}
\put(53,32){\scriptsize{$st_1$}}
\end{picture}
$$
\caption{The dual resolution graph $\Gamma$ of the curve $C$.}
\label{fig1}
\end{figure}

Let us assume that ${\bf 1} \neq st_1$.
Let $\alpha_0={\bf 1}$, $\alpha_1$, \dots, $\alpha_q$ be the dead ends of the
graph $\Gamma$ which do not belong to $\Delta$, and let $\beta_j$
($j=1, \ldots, q$) be the star point of the graph $\Gamma$
which corresponds to the dead end $\alpha_j$, $\beta_1 < \beta_2 < \ldots <
\beta_q$ (see Fig.\ref{fig1}).
Let $S_1$ be the subsemigroup of the semigroup $S_C$ generated by the
multiplicities $\mm^{\alpha_0}$, $\mm^{\alpha_1}$, \dots, $\mm^{\alpha_q}$.
$S_1$ coincides with the semigroup generated by all
the multiplicities $\mm^\sigma$ with $E_\sigma$ from the connected
component of $\overline{D\setminus D^\prime}$ which contains the starting
divisor ${\bf 1}$ and is similar to a subsemigroup
of $\Z_{\ge 0}$ because it is contained in the line $L$ in
$\R^r\supset\Z_{\ge 0}^r$ which goes through the origin and
the point $\mm^{\alpha_0}$.
To describe it more precisely, let $e_i=g.c.d.(m_i^{\alpha_0},
\ldots, m_i^{\alpha_q})$ for $1 \le  i\le r$. Then the set of
integers $\{m_i^{\alpha_0}/e_i$, $m_i^{\alpha_1}/e_i$, \dots,
$m_i^{\alpha_q}/e_i\}$
does not depend on $i$ and is the minimal set of generators of
the semigroup, say $\widetilde{S}$, of an irreducible curve (for
example, of the curve $L_{A}$ with $A\in
{\stackrel{\circ}{E}}_{\beta_q}$). Moreover $\widetilde{m}\in
\widetilde{S}$
if and only if $\widetilde{m}\cdot(e_1, \ldots, e_r)\in S_1$.

As in the case of one branch one has
$\mm^{\beta_j}=(n_j+1)\cdot\mm^{\alpha_j}$, $\mm^{\beta_j}\in
\langle\mm^{\alpha_0}, \ldots, \mm^{\alpha_{j-1}}\rangle$. The
multiplicity $\mm^{st_1}$ belongs to the semigroup $S_1$ as well
(see (5) in \ref{irred}). Let
$\mm^{st_1}=\sum\limits_{j=0}^q  \ell_j\cdot\mm^{\alpha_j}$, where
$\ell_j\le n_j$ for $j=1, \ldots, q$ (such a representation is
unique just as in the case of the semigroup of an irreducible
curve). Let $S_1^\prime$ be the subset of $S_1$ which consists of
elements $\mm$ such that $\mm-\mm^{st_1}\not\in S_1$. $S^\prime_1$
corresponds to the Apery base (see, e.g., \cite{D1}) of
$\widetilde{S}$ with respect to $m_i^{st_1}/e_i$ (this integer
does not depends on $i$). Thus $S^\prime_1$ is a finite set and
the biggest element in it is equal to $\mm^{st_1} +
\sum\limits_{j=1}^{q} n_j \mm^{\alpha_j} - \mm^{\alpha_0}$ (this
follows from the expresion for the conductor in
\ref{irred}).
Let $Y_1=\sum\limits_{\mm\in
S_1^\prime}\bullet\cdot\tt^\mm$. If ${\bf 1} = st_1$ we simply put
$Y_1 = \bullet$.

\smallskip
Let $\widetilde Y =
Y^\prime\times Y_1\times\prod\limits_{\delta\in\Delta} Y_\delta$.

\begin{proposition}\label{prop3}
$\chi{(\widetilde{Y})} = \chi(Y)$.
\end{proposition}

\begin{proof}
One has
$\chi(Y)=\prod\limits_{\sigma\in\Gamma}
\left(1-\tt^{\mm^\sigma}\right)^{-\chi({\stackrel{\circ}{E}}_\sigma)}$.
Since for all $\sigma$ except those from $\Delta^\prime$,
$\Delta$, $\{\alpha_i\}$, and $\{\beta_i\}$,
$\chi({\stackrel{\circ}{E}}_\sigma)=0$,
$$\chi(Y)=\left(\prod\limits_{\sigma\in\Delta^\prime}
\left(1-\tt^{\mm^\sigma}\right)^{-\chi({\stackrel{\circ}{E}}_\sigma)}\right)\times\left(
\prod\limits_{\delta\in\Delta}\left(1-\tt^{\mm^\delta}\right)^{-1}\right)
\times
\frac{\prod\limits_{i=1}^q(1-\tt^{\mm^{\beta_i}})}
{\prod\limits_{i=0}^q(1-\tt^{\mm^{\alpha_i}})}=
$$

$$
=\left(\prod\limits_{\sigma\in\Delta^\prime}
\left(1-\tt^{\mm^\sigma}\right)^{-\chi({\widetilde{E}}_\sigma)}\right)
\times\left(
\prod\limits_{\delta\in\Delta}\frac{(1-\tt^{\mm^{st_\delta}})}
{(1-\tt^{\mm^\delta})}\right)\times
\frac{(1-\tt^{\mm^{st_{\bf 1}}})\prod\limits_{i=1}^q(1-\tt^{\mm^{\beta_i}})}
{\prod\limits_{i=0}^q(1-\tt^{\mm^{\alpha_i}})}.
$$
The first factor coincides with $\chi(Y^\prime)$. Now the statement
follows from the facts that
$$
\chi(Y_\delta)=\frac{1-\tt^{\mm^{st_\delta}}}{1-\tt^{\mm^{\delta}}}, \qquad
\chi(Y_1)=\frac{(1-\tt^{\mm^{st_1}})\prod\limits_{i=1}^q(1-\tt^{\mm^{\beta_i}})}
{\prod\limits_{i=0}^q(1-\tt^{\mm^{\alpha_i}})}.
$$
The first equation is obvious. To prove the second one let us
notice that
$$
\sum\limits_{\vv\in S_1}\tt^\vv=\left(\sum\limits_{\vv\in
S^{\prime}_1}\tt^\vv\right)\cdot\left(1+\tt^{\mm^{st_1}}+\tt^{2\mm^{st_1}}+\ldots\right)
$$
(since each element $s\in S_1$ in a unique way is represented in
the form $\ell\cdot\mm^{st_1}+s^\prime$ with $s^\prime\in S^{\prime}_1$,
$\ell\ge 0$) and
$$
\sum\limits_{\vv\in
S_1}\tt^\vv=\frac{\prod\limits_{i=1}^q(1-\tt^{\mm^{\beta_i}})}
{\prod\limits_{i=0}^q(1-\tt^{\mm^{\alpha_i}})}
$$
(it  follows from the fact that $\mm^{\beta_i}$ are
multiples of
$\mm^{\alpha_i}$:
$\mm^{\beta_i}=(n_i+1)\mm^{\alpha_i}$, $1\le i\le q$, and each element
$s\in S_1$ in a unique way can be represented in the form
$k_0\cdot\mm^{\alpha_0}+\sum\limits_{i=1}^q
k_i\cdot\mm^{\alpha_i}$ with $k_0\ge 0$
and $0\le k_i \le n_i$ for $1\le i\le q$; see, e.g.,
\cite{CDG2}).
\end{proof}

\medskip
There exists a map $\widetilde{\Pi}:
\widetilde{Y}\to \P \widehat{S}_{C}$ such that $Im\,\widetilde{\Pi}
=Im\,\Pi$.
To define it, one can say that $\widetilde Y$ is a subset of the
graded semigroup
$$\widetilde Y^*=Y^\prime\times\left(\sum\limits_{\mm\in S_1}
\bullet\cdot\tt^\mm\right)\times
\prod\limits_{\delta\in\Delta}\left(
\sum\limits_{k=0}^{\infty}\bullet\cdot\tt^{k\mm^\delta}\right)
$$
(each factor of $\widetilde Y^*$ is a graded semigroup)
and the map $\widetilde\Pi$ is the restriction of a graded semigroup
homomorphism $\widetilde Y^*\to\P\widehat S_C$. Because of that it
should be defined for points of
$\bigcup\limits_{\sigma\in\Delta^\prime}\widetilde E_\sigma$ and also for
"monomials" of the form $\bullet\cdot\tt^{\mm^\delta}$ for
$\delta\in\Delta$ and $\bullet\cdot\tt^{\mm^{\alpha_i}}$ for
$i=0, 1, \ldots, q$.
For a point $A$ of ${\stackrel{\circ}{E}}_\sigma$,
$\sigma\in\Delta^\prime$, (or rather for the monomial
$[A]\cdot\tt^{\mm^\sigma}$) $\widetilde\Pi$ coincides with $\Pi$.

A point of $\widetilde E_\sigma\setminus{\stackrel{\circ}{E}}_\sigma$,
$\sigma\in\Delta^\prime$, corresponds either to a dead
end $\delta\in\Delta$ (and in this case $\sigma=st_\delta$) or to the
initial divisor ${\bf 1}$ (in this case $\sigma=st_{1}$). In the
first case one puts $\widetilde\Pi([A]\cdot\tt^{\mm^\sigma})=
(n_\delta+1)\cdot\Pi([A_\delta]\cdot\tt^{\mm^\delta})$ for any
point $A_\delta\in{\stackrel{\circ}{E}}_\delta$; in the
second case one puts $\widetilde\Pi([A]\cdot\tt^{\mm^\sigma})=
\sum\limits_{i=0}^q\ell_i\cdot\Pi([A_{\alpha_i}]\cdot\tt^{\mm^{\alpha_i}})$
for any points $A_{\alpha_i}\in{\stackrel{\circ}{E}}_{\alpha_i}$
(see the definitions of $n_\delta$ and $\ell_i$ above). One puts
$\widetilde\Pi(\bullet\cdot\tt^{\mm^\delta})=
\Pi([A_\delta]\cdot\tt^{\mm^\delta})$ for any point
$A_\delta\in{\stackrel{\circ}{E}}_\delta$, $\delta\in\Delta$,
$\widetilde\Pi(\bullet\cdot\tt^{\mm^{\alpha_i}})=
\Pi([A_{\alpha_i}]\cdot\tt^{\mm^{\alpha_i}})$ for any point
$A_{\alpha_i}\in{\stackrel{\circ}{E}}_{\alpha_i}$, $i=0, 1,
\ldots, q$
(one can easily see that the result does not depend on the choice of the
points $A_\delta$, $A_{\alpha_i}$ in these cases).

\medskip
It is not difficult to see that $Im\,\widetilde{\Pi}
=Im\,\Pi$.

Now the Theorem~\ref{theo4} follows from the following

\begin{proposition}\label{prop4}
$\chi(\widetilde{Y})=\chi(\mbox{Im }\widetilde{\Pi})$.
\end{proposition}

\begin{remark}
Before we prove Proposition~\ref{prop4},
we explain why and in which
sense one can deal with points of $\widetilde E_\sigma
\setminus{\stackrel{\circ}{E}}_\sigma$, $\sigma\in\Delta^\prime$,
just in the same way as with
other points of $\widetilde E_\sigma$. Let $A\in\widetilde E_\sigma
\setminus{\stackrel{\circ}{E}}_\sigma$. The point $A$ corresponds either
to a dead end $\delta\in\Delta$ (in this case $\sigma=st_\delta$)
or to the starting divisor ${\bf 1}$ (in this case $\sigma=st_{1}$).
Let $A^\prime\in{\stackrel{\circ}{E}}_\sigma$, and let $g$ and $g^\prime$
be functions $(\C^2, 0)\to(\C, 0)$ such that
$\widetilde\Pi([A]\cdot\tt^{\mm^\sigma})=(\vv(g), \aa(g))$,
$\widetilde\Pi([A^\prime]\cdot\tt^{\mm^\sigma})=(\vv(g^\prime),
\aa(g^\prime))$.
Let us recall that $g^\prime=g_{A^\prime}$, $g=
(g_{A_\delta})^{n_\delta}$ with
$A_\delta\in{\stackrel{\circ}{E}}_\delta$ if the
point $A$ corresponds to the dead end $\delta\in\Delta$, $g=
(g_{A_{\alpha_0}})^{\ell_0}\cdot(g_{A_{\alpha_1}})^{\ell_1}\cdot\ldots\cdot
(g_{A_{\alpha_q}})^{\ell_q}$ with
$A_{\alpha_i}\in{\stackrel{\circ}{E}}_{\alpha_i}$ if the
point $A$ corresponds to the starting divisor ${\bf 1}$  (see the
notations above). Then $\vv(g)=\vv(g^\prime)$. To compare $\aa(g)$
and $\aa(g^\prime)$ one can look (and we shall
regularly do it below)
at the ratio $\psi={\widetilde g^\prime}/{\widetilde g}$, where
$\widetilde g = g\circ\pi$ and $\widetilde g^\prime =
g^\prime\circ\pi$ are the liftings of the functions $g$ and $g^\prime$
to the space $X$ of the resolution. The main (or rather the only) property
of the function $\psi$ which will be used below is the following one.
The restriction $\psi_{\vert E_\sigma}$ of the function $\psi$ to
the component $E_\sigma$ of the exceptional divisor is a meromorphic
function (in fact a ratio of two linear functions) with one pole at
the point $A$ and one zero at the point $A^\prime$. It has no zeroes
or poles on all other components $E_{\sigma^\prime}$ from $D^\prime$.
Therefore it is constant (and different from zero or infinity) on each
connected component of
$D^\prime\setminus E_\sigma$ and its value on such a component
coincides with the value of the function $\psi$ at the corresponding
essential point of the component $E_\sigma$.
\end{remark}

{\bf Proof of Proposition~\ref{prop4}.}
The proof consists in analyzing places where $\widetilde{\Pi}$ is not
immersive. At each such place we explicitly show a part of
$\widetilde{Y}$ which has the Euler characteristic equal to zero and
which can be removed without changing the image.

\smallskip
Let us write $\widetilde Y=\sum_{\kk}Y_\kk$, where $\kk$
is the multi-index
$$
\kk = \left\{
\{k_{\sigma}\}_{\sigma\in \Delta'}, \mm,
\{k_{\delta}\}_{\delta\in \Delta} \right\}
$$
with
$k_{\sigma}\ge 0$ for each  $\sigma\in \Delta'$, $\mm\in
S_1^\prime$,
and $0\le k_{\delta}\le n_{\delta}$ for each $\delta\in \Delta$.
One has $ Y_\kk=
(\prod\limits_{\sigma\in\Delta^\prime}{S^{k_\sigma}\widetilde
E_\sigma})\cdot \tt^{\vv(\kk)}$,
where
$\vv(\kk)=\sum k_\sigma\mm^\sigma+\mm+\sum k_\delta\mm^\delta$
($Y_\kk$ are "connected components" of $\widetilde Y$).

Suppose that there exist two different elements $y_1$ and $y_2$ from
$\widetilde Y$ such that $\widetilde\Pi(y_1)=\widetilde\Pi(y_2)$.
Let  $y_j\in Y_{\kk_j}$ ($j=1,\, 2$), and let the element $y_j$
be represented
by the sets $K_{\sigma,j}\subset\widetilde E_\sigma$ for
$\sigma\in\Delta^\prime$ ($\# K_{\sigma,j}=k_{\sigma,j}$),
i.e., $y_j = (\prod\limits_{\sigma\in\Delta^\prime}([K_{\sigma,j}]
\cdot\tt^{k_{\sigma,j}\mm^\sigma})\times x_j$,
where $x_j\in Y_1\times\prod\limits_{\delta\in\Delta}Y_{\delta}$.

\begin{lemma}\label{lemma3}
There exists a component $E_\sigma$ of ${D^\prime}$ such that:
\newline
1) for all $\sigma^*\in\Delta^\prime$, $\sigma^*>\sigma$, one has
$K_{\sigma^*,1} = K_{\sigma^*,2}$; \newline
2) for all $\delta\in\Delta$, $\delta>\sigma$, one has
$k_{\delta,1} = k_{\delta,2}$; \newline
3) $K_{\sigma,1} \ne K_{\sigma,2}$ and either $k_{\sigma,1}$ or
$k_{\sigma,2}$ is $\ge s_\sigma-1$.
\end{lemma}

\begin{proof}
First let us show that there exists a component
$E_\sigma\subset D^\prime$ such that $K_{\sigma,1} \ne K_{\sigma,2}$.
If $K_{\sigma,1} = K_{\sigma,2}$ for all $E_\sigma\subset D^\prime$,
then $\widetilde\Pi(x_1)=\widetilde\Pi(x_2)$ ($x_j\in
Y_1\times\prod\limits_{\delta\in\Delta}Y_\delta$; see above). Let
$x_j=x_j^{(1)}\times\prod\limits_{\delta\in\Delta}x_j^{(\delta)}$, where
$x_j^{(1)}\in Y_1$, $x_j^{(\delta)}\in Y_\delta$, and let us
suppose that $x_1^{(\delta_0)}\neq x_2^{(\delta_0)}$ for
$\delta_0\in\Delta$, but $x_1^{(\delta)} = x_2^{(\delta)}$ for
all dead ends $\delta$ such that $st_{\delta_0} < \delta$.
Without loss of generality one
can suppose that $x_j^{(\delta)}=\bullet$ for those $\delta$
($j=1,2$) and $x_2^{(\delta_0)}=\bullet$,
$x_1^{(\delta_0)}=\bullet\cdot \tt^{k\mm^{\delta_0}}$ with $0<k\le
n_{\delta_0}$.

Let $C_i$ be a component of the curve $C$ such that
$st_{\delta_0} < \widetilde C_i$. In this case
$E_{\delta_0}$ and $E_{st_{\delta_0}}$ appear in the minimal
embedded resolution of the (irreducible) curve $C_i$ and as a
consequence
$m_i^{\delta_0}$ (the $i$-th component of $\mm^{\delta_0}$)
belongs to the minimal set of generators
$\{\bar\beta_0, \bar\beta_1, \ldots, \bar\beta_e\}$ of the
semigroup $S_{C_i}$ ($\subset\Z_{\ge 0}$) of the
curve $C_i$: $m_i^{\delta_0}=\bar\beta_p$
(with the notations of  \ref{irred},
in the dual graph of the minimal resolution of the curve $C_i$ one has
$\delta_0 = \alpha_p$ and
$st_{\delta_0}=\beta_p$).
Moreover, for any $\delta\in \Delta$ with
$\delta\not > st_{\delta_0}$   (i.e., either $\delta <
st_{\delta_0}$ or $\delta$ and $st_{\delta_0}$ are not
comparable) the function
$g_{A_{\delta}}$ satisfies the hypothesis of the point (5) in
\ref{irred} and so
the corresponding value
$v_i(g_{A_{\delta}})\in S_{C_i}$ belongs to
$\langle
\bar\beta_0, \bar\beta_1, \ldots, \bar\beta_{p-1}
\rangle$.
As a consequence one has that
$v_i(\widetilde\Pi(x_1))=k\bar\beta_p+v^\prime$ where $0<k\le
n_{\delta_0}$, $v^\prime\in\langle\bar\beta_0, \bar\beta_1,
\ldots, \bar\beta_{p-1}\rangle$,
$v_i(\widetilde\Pi(x_2))\in\langle\bar\beta_0, \bar\beta_1,
\ldots, \bar\beta_{p-1}\rangle$
(here $v_i$ is a coordinate in
$\P\widehat S_C \subset\Z_{\ge 0}^r\times\P((\C^*)^r)$\,).
This contradicts the uniqueness of the representation in the semigroup
$S_{C_i}$ (see properties (2) and (4) of
\ref{irred}).
This proves the statement in the discussed case.
The same arguments (applied to the semigroup $S_1$) work in the
case if $x_1^{(\delta)}=x_2^{(\delta)}$ for all $\delta\in\Delta$.

Let $\sigma$ be a maximal element in the set of vertices
from $\Delta^\prime$ with $K_{\sigma,1} \ne K_{\sigma,2}$, i.e.,
$K_{\sigma^\prime,1} = K_{\sigma^\prime,2}=K_{\sigma^\prime}$
for all  $\sigma^\prime>\sigma$, $\sigma^\prime\in\Delta^\prime$.
Just the previous arguments show that, for all dead ends $\delta$
with $st_\delta\ge\sigma$, one has $k_{\delta, 1}=k_{\delta, 2}$.
Let
$$
y^\prime_j =\left(
\prod\limits_{\sigma^\prime\in\Delta^\prime,\,\sigma^\prime\not>
\sigma}
[K_{\sigma^\prime,j}]\cdot\tt^{k_{\sigma^\prime,j}\mm^{\sigma^\prime}}
\right)
\times
x^\prime_j,
$$
$j=1,2$, where $\sigma^\prime\not> \sigma$ means that either
$\sigma^\prime\le\sigma$ or $\sigma$ and $\sigma^\prime$ are not comparable,
$x^\prime_j = \prod\limits_{\delta\in\Delta,\,\delta\not>
\sigma}(\bullet\cdot\tt^{k_{\delta,j}\mm^{\delta}})\times
(\bullet\cdot\tt^{\mm_j})$, $\mm_j\in S_1^\prime$ (i.e.,
$x^\prime_j$ is obtained from $x_j$ by dropping all factors
$\tt^{k_{\delta,j}\mm^{\delta}}$ with $\delta > \sigma$).
One has $y^\prime_1 \neq y^\prime_2$,
$\widetilde{\Pi}(y^\prime_1)=\widetilde{\Pi}(y^\prime_2)$
(the last equation follows from the fact that multiplication (the
semigroup operation) by any element of $\Z_{}^r\times(\C^*)^r$ is injective).

Let $k_j = k_{\sigma,j} = \# K_{\sigma,j}$, $j = 1,2$. Without
loss of generality  one can suppose that $k_1\ge k_2$.
Let $Q_0$, $Q_1$, \dots $Q_{s-1}$ ($s=s_\sigma$) be essential points
on the component $E_\sigma$. If $\sigma\ne st_1$, we suppose that
the point $Q_0$ corresponds to the connected component of
$D\setminus {\stackrel{\circ}{E}}_\sigma$ which contains the
starting divisor~${\bf 1}$.

Let us fix an affine coordinate on the projective line $E_{\sigma}$ in
such a way that the essential point $Q_0$ of $E_{\sigma}$ (corresponding
to the starting divisor ${\bf 1}$ if $\sigma\ne st_1$) is the infinite one.

Let $g_1$ and $g_2$ ($g_j: (\C^2,0)\to  (\C,0)$) be functions,
corresponding to $y^\prime_1$ and $y^\prime_2$, let
$\widetilde{g}_j=g_j\circ \pi$ be the lifting of the
function $g_j$ to the space ${X}$ of the resolution, and let
$\psi = {\widetilde{g}_1}/{\widetilde{g}_2}$. The function $\psi$
is a meromorphic function on $E_{\sigma}$, and is a regular
nonzero function on
$\bigcup\limits_{\sigma^\prime>\sigma} E_{\sigma^\prime}$. Therefore it is
constant on each connected component of
$\bigcup\limits_{\sigma^\prime>\sigma} E_{\sigma^\prime}$. The value of
$\psi$ on such a component is equal to $\psi(Q_\ell)$, where
$Q_\ell$ is the essential point of $E_{\sigma}$
corresponding to the component.

Since $\aa(g_1) = c\cdot\aa(g_2)$, $\psi(Q_\ell) = c$ for a constant
$c \neq 0$, $\ell=1, \ldots, s-1$. On $E_{\sigma}$ the function
$\psi$ has the form $\psi = c^\prime\cdot\frac{p_1(z)}{p_2(z)}$, where
$p_j(z)= \prod\limits_{k=1}^{k_j} (z-z_k^{(j)})$,
$K_{\sigma,j}= \{z_k^{(j)}\}$, $c^\prime \neq 0$.

The polynomial $p(z) = c^\prime\cdot p_1(z)-c\cdot p_2(z)$ vanishes
at all the points $Q_\ell$ ($\ell=1, \ldots, s-1$). Since
$K_{\sigma,1} \neq K_{\sigma,2}$, $p(z) \neq 0$. One has
$\deg p(z)\le k_1$ and $p(z)$ has (at least) $s-1$
zeroes. Therefore $k_1\ge s-1$.
\end{proof}

\begin{definition}
We shall say that a vertex $\sigma\in\Delta^\prime$ is {\bf a cut}
of a multi-index $\kk=\{k_\sigma, \mm, k_\delta\}$
if $k_\sigma\ge s_\sigma-1$.
\end{definition}

Lemma~\ref{lemma3} implies that, if
there exist $y_1\in Y_{\kk_1}$, $y_2\in Y_{\kk_2}$, $y_1\ne y_2$,
such that $\widetilde\Pi(y_1)=\widetilde\Pi(y_2)$, then either
$\kk_1$ or $\kk_2$ has a cut.

\begin{remark}
Assume that a multi-index $\kk$ has a cut at $\sigma\in \Delta'$.
Then the Euler characteristic of the component $Y_{\kk}$ of the
space $\w{Y}$ is equal to zero and thus it makes no contribution
to $\chi(\w{Y})$.
By lemma \ref{lemma3}, if $\w{\Pi}(Y_{\kk_1})\cap
\w{\Pi}(Y_{\kk_2}) \neq \emptyset$, then either $\kk_1$ or
$\kk_2$ has a cut at some place.
The idea is that one can
remove the component with a cut (at least part of it) without
changing the image. The rest of the proof explains
how this can be made. The most technical part consists in showing
some fine properties of cuts and their relative distribution in
the dual graph (Lemmas \ref{lemma5} and \ref{lemma6}, Lemma
\ref{lemma4} is a technical tool to simplify the proof of the
others).
\end{remark}

\subsection{Distribution of cuts.}\label{s-cuts}
Let
$$
y_j = \prod\limits_{\sigma\in\Delta^\prime}
\left([K_{\sigma,j}]\cdot\tt^{k_{\sigma,j}\mm^\sigma}\right)\times
(\bullet\cdot\tt^{\mm_j})\times\prod\limits_{\delta\in \Delta}
(\bullet\cdot\tt^{k_{\delta,j}\mm^\delta})
$$
($j=1, 2$, $\# K_{\sigma,j}=k_{\sigma,j}$, $\mm_j\in S_1^\prime$,
$0\le k_{\delta,j}\le n_\delta$) and let
$$
g_j =
\prod\limits_{\sigma\in\Delta^\prime}\varphi_j^{(\sigma)}
\cdot \varphi_j^{(1)}\cdot
\prod\limits_{\delta\in \Delta} \varphi_j^{(\delta)}
$$
be the corresponding function from ${\cal O}_{\C^2, 0}$:
$\widetilde\Pi(y_j)=(\vv(g_j), \aa(g_j))$.

For a vertex $\sigma\in \Delta^\prime$, let $Q_0$, $Q_1$, \dots
$Q_{s-1}$ ($s=s_{\sigma}$) be the essential points of
$E_{\sigma}$. Let $\Gamma(\sigma,\ell)$ be the subgraph of
$\Gamma$ corresponding to the  connected component of
$D\setminus {\stackrel{\circ}{E}}_\sigma$ which contains the point $Q_\ell$. If
$\sigma\ne st_1$, we assume that ${st_1}\in \Gamma(\sigma,0)$.
In what follows we shall use the following notations:
$$
\begin{array}{rl}
\vv({g_1}/{g_2}; \ge \sigma) &=
\sum\limits_{\tau\in \Delta^\prime\cup \Delta, \tau\ge \sigma }
(\vv(\varphi_1^{(\tau)})-
\vv(\varphi_2^{(\tau)})), \\
\vv({g_1}/{g_2}; > \sigma) &=
\sum\limits_{\tau\in \Delta^\prime\cup \Delta, \tau>\sigma}
(\vv(\varphi_1^{(\tau)})-
\vv(\varphi_2^{(\tau)})), \\
\vv({g_1}/{g_2}; \sigma,\ell) &=
\sum\limits_{\tau\in \Gamma(\sigma,\ell)}
(\vv(\varphi_1^{(\tau)})-
\vv(\varphi_2^{(\tau)})) \\
\end{array}
$$
(here $\vv(\varphi_1^{(\tau)}) - \vv(\varphi_2^{(\tau)})=
(k_{\tau,1} - k_{\tau,2})\cdot\mm^\tau$, $\tau\in \Delta^\prime\cup \Delta$).
For $\sigma>st_1$ (respectively for $\sigma=st_1$),
$\vv(g_1/g_2; > \sigma)$ is equal to $\sum\limits_{\ell=1}^{s-1}
\vv(g_1/g_2;\sigma,\ell)$ (respectively to
$\sum\limits_{\ell=0}^{s-1} \vv(g_1/g_2;\sigma,\ell)$\,) plus the
contribution, $\vv(\varphi^{\delta}_1) -
\vv(\varphi^{\delta}_2)$, of the dead end $\delta$ such that
$st_\delta=\sigma$ if such a dead end exists.

\medskip
Let $\sigma\in \Delta^\prime$. We define a ``modified"
multiplicity, $\w{\mm}^{\sigma}=(\w{m}^{\sigma}_1, \ldots ,
\w{m}^{\sigma}_r)$, in the following way. If $\sigma$ is a
star vertex and there exists a dead end $\delta\in \Delta$ such
that
$st_{\delta}\ge \sigma$ and $\mm^{\sigma}>\mm^{\delta}$, we put
$\w{\mm}^{\sigma} = \mm^{\delta}$; otherwise $\w{\mm}^{\sigma}
= \mm^{\sigma}$ (a dead end $\delta$ with the described
properties, if it exists, is unique). The dead end
$\delta\in\Delta$ with $st_\delta\ge\sigma$, $\mm^\sigma >
\mm^\delta$ (if it exists) will be denoted by $\delta(\sigma)$.
Such dead end $\delta$ exists if and only if, in the process of
resolution by blow-ups of points, the component $E_\sigma$ of
the exceptional divisor is created after the component $E_{\delta}$
(i.e., by blowing-up points of $E_{\delta}$)
and $\sigma\le  st_{\delta}$.

The reason for this construction is that the natural order
in the dual graph (that is the order we use) does not coincides
with the order in which the successive divisors are created.
As a consequence the multiplicity map $\alpha \mapsto \mm^\alpha$
is not an increasing one.

\begin{remarks}
{\bf 1.} If $\sigma=st_{\delta}$ for $\delta\in\Delta$, then
$\delta=\delta(\sigma)$. However $\delta(\sigma)$ could also exist
in the case when $\sigma \neq st_{\delta}$ for any $\delta\in
\Delta$.
Moreover, if $\delta(\sigma)$ exists for some
$\sigma$ then $\delta(\sigma) = \delta(\sigma^\prime)$ for all star
vertices $\sigma^\prime\in \Delta^\prime$ with $\sigma\le \sigma^\prime\le
st_{\delta}$. Therefore one and the same dead end $\delta$ could occur as
$\delta(\sigma)$ for several star vertices $\sigma\in
\Delta^\prime$.

{\bf 2.} Let $\sigma$ be a star vertex. Then for any $\tau\in
\Gamma$ with $\tau>\sigma$ one has that $\mm^{\tau}\ge
\w{\mm}^{\sigma}$ and $\mm^{\tau} = \w{\mm}^{\sigma}$ if and only
if $\tau=\delta(\sigma)$. Therefore $\w{\mm}^{\tau}\ge
\w{\mm}^{\sigma}$ for any $\tau>\sigma$ and the equality holds if
and only if either $\tau$ is a star vertex and
$\delta(\tau)=\delta(\sigma)$ or $\tau$ is a dead end and
$\tau=\delta(\sigma)$.
\end{remarks}

\begin{lemma}\label{lemma4}
Let $\sigma\in \Delta^\prime$, $\sigma>st_1$, be a vertex such
that $\vv(g_1/g_2; \ge \sigma)\ge\w{\mm}^{\sigma}$. Let $\sigma^\prime$
be the previous separation
vertex of $\Gamma$, i.e., either $\sigma^\prime=st_1$ or
$s_{\sigma^\prime}>2$, $\sigma^\prime<\sigma$ and
$s_{\sigma''} = 2$ for any $\sigma''\in \Delta^\prime$ with
$\sigma^\prime<\sigma''<\sigma$. Let $\ell_0$ be such that
$\sigma\in \Gamma(\sigma^\prime, \ell_0)$.
Assume that there are no cuts of $\kk_2$ between $\sigma^\prime$ and
$\sigma$. Then: \newline
--- if there exists $\delta(\sigma^\prime)$ and
it belongs to $\Gamma(\sigma^\prime, \ell_0)$,
then $\vv(g_1/g_2; \sigma^\prime, \ell_0)\ge \w{\mm}^{\sigma^\prime}$;
\newline
--- otherwise $\vv(g_1/g_2; \sigma^\prime, \ell_0)\ge \mm^{\sigma^\prime}$.

Moreover, let $i,j\in \{1,\ldots,r\} $ be such that $\sigma <
\w{C}_i$ and $\sigma\not<\w{C}_j$.
Assume that $v_j(g_1/g_2; \ge \sigma)\cdot m_i^{\sigma} \le
v_i(g_1/g_2; \ge \sigma)\cdot m_j^{\sigma}$.  Then
$v_j(g_1/g_2; \sigma^\prime,\ell_0)\cdot m_i^{\sigma^\prime} <
v_i(g_1/g_2; \sigma^\prime,\ell_0)\cdot m_j^{\sigma^\prime}$.
\end{lemma}

\begin{proof}
Let us prove the first statement. Let $\{\delta_1, \delta_2,
\ldots , \delta_p\}$ be (all) dead ends such that $\sigma^\prime
< st_{\delta_1}<\ldots < st_{\delta_p} < \sigma$. Suppose first that
$p=0$. If there exists $\delta(\sigma^\prime)$ and it belongs to
$\Gamma(\sigma^\prime, \ell_0)$, then either $\delta(\sigma) =
\delta(\sigma^\prime)$ (if $\sigma$ is a star vertex) and
$\w{\mm}^{\sigma^\prime} =\w{\mm}^{\sigma}$ or $\w{\mm}^{\sigma} =
\mm^{\sigma} > \w{\mm}^{\sigma^\prime}$. In both cases the
statement is obvious. If $\delta(\sigma^\prime)$ does not belong
to $\Gamma(\sigma^\prime, \ell_0)$ (or does not exist), then
$\w{\mm}^{\sigma^\prime} < \w{\mm}^{\sigma}$ and the statement is
obvious as well.

Suppose that $p>0$. Since $\delta(st_{\delta_p})= \delta_p\not>
\sigma$, $\w{\mm}^{\sigma} > \mm^{st_{\delta_p}} =
(n_{\delta}+1)\cdot\mm^{\delta_p}$. Taking into account that
$\vv(\varphi_2^{(\delta_p)}) \le n_{\delta_p}\mm^{\delta_p}$, one
has
$$
\begin{array}{rl}
\vv({g_1}/{g_2}; \ge st_{\delta_p}) &\ge
\vv({g_1}/{g_2}; \ge\sigma) -
n_{\delta_p}\mm^{\delta_p} > \\
\ & > (n_{\delta_p}+1)\mm^{\delta_p} -
n_{\delta_p}\mm^{\delta_p} = \mm^{\delta_p}
=\w{\mm}^{st_{\delta_p}}.
\end{array}
$$
Repeating the same arguments for the star points $st_{\delta_{p-1}}$,
$st_{\delta_{p-2}}$, \dots, one proves the statement.

Now we prove the second statement. Let $\sigma^*$ be the
separation vertex of the branches $C_i$ and $C_j$ (one has
$\sigma^*\le \sigma^\prime$). For $\tau\in \Gamma$, let
$h_\tau = m_j^{\tau}/m_i^{\tau}$. If $S_{ij}\subseteq \Z_{\ge
0}^{2}$ is the semigroup of the curve $C_i\cup C_j$ (it
coincides with the projection of the semigroup $S=S_C$ to the
$(v_i, v_j)$--plane), then $h_\tau$ is just the slope of the
line which goes through the origin and the point
$(m_i^\tau,m_j^\tau)\in S_{ij}$. It is known
(see \cite{CDG0}, proof of Theorem 2)
that the slopes $h_\tau$ are constant for
$\tau\in
\Delta^\prime$ with ${\bf 1}\le \tau\le \sigma^*$; decrease strictly
for $\tau$ such that $\sigma^* \le \tau\le \w{C}_i$ (i.e.,
$\sigma^*\le \tau< \tau^\prime\le \w{C}_i$ if and only if
$h_\tau > h_{\tau^\prime}$) and increase strictly for $\tau$ such
that $\sigma^* \le \tau\le \w{C}_j$. Moreover $h_\tau$ is
constant on each tail. If $p=0$, the statement is obvious because
$v_u(g_1/g_2; \sigma^\prime, \ell_0)=v_u(g_1/g_2; \ge\sigma)$
for $u=i, j$ and $m_j^{\sigma^\prime}/m_i^{\sigma^\prime}=
h_{\sigma^\prime} > h_{\sigma} = m_j^{\sigma}/m_i^{\sigma}$.

If $p > 0$, we have the inequalities
$$
\frac{v_j(g_1/g_2;\ge \sigma)}{v_i(g_1/g_2;\ge \sigma)} \le
\frac{m_j^{\sigma}}{m_i^{\sigma}} <
\frac{m_j^{\delta_p}}{m_i^{\delta_p}} < \cdots <
\frac{m_j^{\delta_1}}{m_i^{\delta_1}} <
\frac{m_j^{\sigma^\prime}}{m_i^{\sigma^\prime}}\, ,
$$
$$
\ww^\sigma:=(v_i(g_1/g_2; \ge\sigma),v_j(g_1/g_2; \ge \sigma))
\ge
(\w{m}_i^{\sigma},\w{m}_j^{\sigma}) >
(n_{\delta_p}+1)
(m_i^{\delta_p},m_j^{\delta_p})\;.
$$
Let $\ww^{st_{\delta_p}} = (v_i(g_1/g_2; \ge
st_{\delta_p}),v_j(g_1/g_2; \ge st_{\delta_p}))$. One has
$\ww^{st_{\delta_p}} = \ww^{\sigma}-\mm$, where $\mm =
(k_{\delta_p,2}-k_{\delta_p,1})\cdot(m_i^{\delta_p},m_j^{\delta_p})$.
Since $\vert k_{\delta_p,2}-k_{\delta_p,1}\vert \le
n_{\delta_p}$, the slope of the vector $\ww^{st_{\delta_p}}$
is less than $h_{\delta_p}$ and thus less than $h_{\delta_{p-1}}$
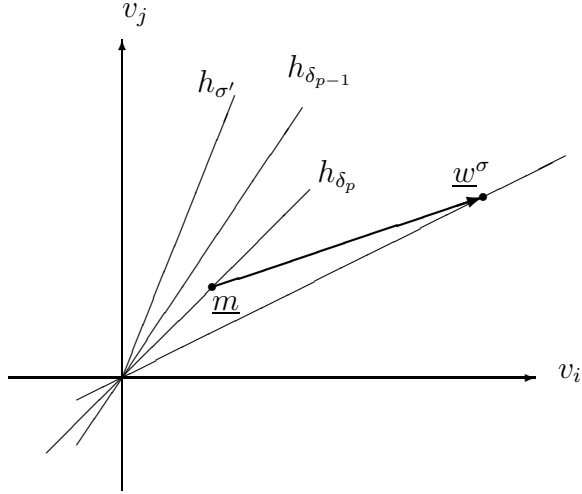
\begin{figure}
$$
\unitlength=1mm
\begin{picture}(80.00,70.00)(-5,-10)
\thinlines
\put(-15,0){\vector(1,0){70}}
\put(0,-15){\vector(0,1){60}}
\put(0,0){\line(2,5){15}}
\put(-6,-9){\line(2,3){30}}
\put(-10,-10){\line(1,1){35}}
\put(-6,-3){\line(2,1){65}}
\put(12,12){\circle*{1}}
\put(48,24){\circle*{1}}
\thicklines
\put(12,12){\vector(3,1){36}}
\put(58,0){$v_i$}
\put(0,48){$v_j$}
\put(10,38){$h_{\sigma^\prime}$}
\put(22,40){$h_{\delta_{p-1}}$}
\put(26,26){$h_{\delta_p}$}
\put(12,9){$\mm$}
\put(44,26){$\ww^{\sigma}$}
\end{picture}
$$
\caption{$(v_i, v_j)$--plane.}
\label{fig2}
\end{figure}
(see Fig.\ref{fig2}). Moreover
$\ww^{st_{\delta_p}}>(m_i^{\delta_p},m_j^{\delta_p})
 >(n_{\delta_{p-1}}+1)(m_i^{\delta_{p-1}},m_j^{\delta_{p-1}})$.
Thus one has:
$$
\frac{v_j(g_1/g_2;\ge st_{\delta_p})}{v_i(g_1/g_2;\ge
st_{\delta_p})} <
\frac{m_j^{\delta_{p-1}}}{m_i^{\delta_{p-1}}}
\quad  {\rm and } \quad
\ww^{st_{\delta_p}} >
(n_{\delta_{p-1}}+1)
(m_i^{\delta_{p-1}},m_j^{\delta_{p-1}})\;  .
$$
Repeating the same arguments for the star points
$st_{\delta_{p-1}}$, $st_{\delta_{p-2}}$, \dots, one proves the
statement.
\end{proof}

\begin{remarks}
Let $\tau\in \Delta^\prime$ be such that $\tau$ is not a cut of
$\kk_2$. If there exists $\delta(\tau)$, let $\ell^*\in \{0, 1,
\ldots, s_{\tau}-1\}$ be such that $\delta(\tau)\in
\Gamma(\tau,\ell^*)$. If the conclusions of the Lemma \ref{lemma4}
are valid for $\tau$, $g_1/g_2$ and for any $\ell=1,\ldots
,s_{\tau}-1$ ($\ell=0,1,\ldots ,s_{\tau}-1$ if $\tau=st_{1}$) then
the hypothesis on the Lemma \ref{lemma4} are true for $\tau$ and
$g_1/g_2$ as well. More explicitly:

{\bf 1.}
Suppose that $\tau \neq st_1$ and
$\vv(g_1/g_2; \tau, \ell)\ge \mm^{\tau}$ for
$\ell = 1,\ldots, s_{\tau}-1$, $\ell\neq \ell^*$ (if it exists)
and $\vv(g_1/g_2; \tau, \ell^*)\ge \w{\mm}^{\tau}$.
Then, if $\tau \neq st_{\delta}$ for any dead end $\delta$,
$\vv(g_1/g_2;> \tau) = \sum\limits_{\ell=1}^{s_{\tau}-1}
\vv(g_1/g_2;\tau, \ell)\ge
(s_{\tau}-2)\cdot\mm^{\tau}+\w{\mm}^{\tau}$.
If $\tau=st_{\delta}$
(in this case $\ell^*$ does not exists) one has
$\vv(g_1/g_2;> \tau) \ge
(s_{\tau}-1)\cdot\mm^{\tau}- \vv(\varphi^{(\delta)}_2) \ge
(s_{\tau}-2)\cdot\mm^{\tau}+\w{\mm}^{\tau}$.

Since $\tau$ is not a cut of $\kk_2$, one has that
$\vv(\varphi_2^{(\tau)})\le (s_{\tau}-2)\cdot\mm^{\tau}$.
Thus,
in any case one has that $\vv(g_1/g_2;\ge \tau)\ge
\w{\mm}^{\tau}$.

If $\tau=st_1$, the only difference in the discussions above
is that the index
$\ell=0$ plays the same role as the others and so
one has $\vv(g_1/g_2;\ge st_1) \ge \mm^{st_1} + \w{\mm}^{st_1}$.

{\bf 2.} Let $i,j\in
\{1,\ldots,r\}$ be such that
$\tau < \w{C}_i$, $\tau\not< \w{C}_j$. Suppose that
$v_j(g_1/g_2; \tau,\ell)\cdot m_i^{\tau} <
v_i(g_1/g_2; \tau,\ell)\cdot m_j^{\tau}$
for
$\ell= 1,\ldots, s_{\tau}-1$.
Then $v_j(g_1/g_2; \ge \tau)\cdot m_i^{\tau} <
v_i(g_1/g_2; \ge \tau)\cdot m_j^{\tau}$.
\end{remarks}

\begin{lemma}\label{lemma5}
Let
$\kk_1\neq \kk_2$ be such that there exist $y_1\in Y_{\kk_1}$
and
$y_2\in Y_{\kk_2}$ with
$\w{\Pi}(y_1) = \w{\Pi}(y_2)$. Suppose that $\sigma_2\in \Delta^\prime$
is such that: \newline
1) there are no cuts $\sigma^\prime_2$ of $\kk_2$ with
$\sigma^\prime_2>\sigma_2$; \newline
2) there exists a cut $\sigma_1$ of $\kk_1$ with
$\sigma_1\ge \sigma_2$. \newline
Then on each geodesic from the vertex $\sigma_2$ to a strict transform
$\w{C}_j$ with $\w{C}_j>\sigma_2$ there exists a cut of $\kk_1$.
\end{lemma}

\begin{proof}
We use the induction on the number $q$ of branches $C_j$ such
that $\w{C}_j >\sigma_2$. The statement is obvious for $q=1$. Let
$q>1$, and let $\w{C}_i, \w{C}_j$ ($i\ne j$) be such that
$\sigma_2 < \sigma_1 < \w{C}_i$, $\sigma_2 < \w{C}_j$.

Suppose that there is no cut of $\kk_1$ on the geodesic from
$\sigma_2$ to $\w{C}_j$, and let $\sigma^*$ be the separation
vertex between $C_i$ and $C_j$. Without loss of generality one
can assume that $\sigma_1$ is maximal among  the cuts of $\kk_1$
on the geodesic from $\sigma_2$ to $\w{C}_i$ and that $\sigma_2$
is a separation vertex. If $\sigma^* >\sigma_2$, or if there is
a cut $\sigma^\prime_1$ of $\kk_1$ on the connected component of
$\overline{D\setminus E_{\sigma^*}}$ which intersects $\w{C}_j$
(this cut must be not comparable with $\widetilde C_j$), the statement
follows from the inductive hypothesis (in the first case the
number of branches $C_j$ such that $\w{C}_j>\sigma^*$ is
strictly smaller than $q$; in the second case one can apply the
arguments to a separation point $\sigma^{**}$ such that
$\sigma^{*}<\sigma^{**}<\widetilde C_j$).

Thus one can assume that $\sigma^*=\sigma_2$ and that
neither $\kk_1$ nor $\kk_2$ has a cut on the connected component
$\Gamma(\sigma^*, \ell_j)$ of
$\Gamma-\{\sigma^*\}$ which contains $\w{C}_j$. Let
$\Gamma(\sigma^*, \ell_i)$ be the connected component of
$\Gamma-\{\sigma^*\}$ which contains $\w{C}_i$ and let
$g_u = G^0_u\cdot G^i_u \cdot G^j_u$ ($u=1,2$) where
$$
G^i_u = \prod\limits_{\tau\in \Gamma(\sigma^*, \ell_i)}
\varphi^{(\tau)}_u\; , \quad
G^j_u = \prod\limits_{\tau\in \Gamma(\sigma^*, \ell_j)}
\varphi^{(\tau)}_u\; , \quad
G^0_u = \frac{g_u}{G^i_u\cdot G^j_u}\; .
$$
Let $\ww(g):=(v_i(g), v_j(g))\in\Z_{\ge 0}^2$.

For any $\tau\in \Delta^\prime\cap \Gamma(\sigma^*, \ell_j)$ one has
$K_{\tau,1} = K_{\tau,2}$
(see Lemma~{\ref{lemma3}}). A simple computation (similar
to the one in the proof in
Lemma~\ref{lemma3}) shows that $x_1^{(\delta)}=x_2^{(\delta)}$
for any $\delta\in \Delta$ with
$st_{\delta}\in \Gamma(\sigma^*, \ell_j)$. Therefore
$$
\ww(G^j_1) = \ww(G^j_2).
$$
Since $\sigma_1$ is a maximal cut of $\kk_1$,
$\vv(g_1/g_2; \ge \sigma_1)\ge \mm^{\sigma_1}$.
Since $\Gamma(\sigma^*,\ell_i)$ does not contain cuts of $\kk_2$,
by the inductive hypothesis,  on each geodesic from the vertex
$\sigma_2$ to a strict transform $\w{C}_n$ with
$\w{C}_n >\sigma_2$, $\w{C}_n\in \Gamma(\sigma^*,\ell_i)$, there
exists a cut of $\kk_1$. Thus, using repeatedly
Lemma~\ref{lemma4} and the Remark after it, one can see that
$$
\begin{array}{rl}
\ww(G^i_1/G^i_2) &\ge (\widetilde m_i^{\sigma^*},
\widetilde m_j^{\sigma^*})>(0,0)\,,\\
(v_j(G^i_1)-v_j(G^i_2))\cdot m_i^{\sigma^*} &<
(v_i(G^i_1)-v_i(G^i_2))\cdot m_j^{\sigma^*}\,.  \\
\end{array}
$$
It is known that, for any factor $\varphi$ of $G^0_u$ ($u=1,2$),
one has
$$
v_j(\varphi)\cdot m_i^{\sigma^*} =
v_i(\varphi)\cdot m_j^{\sigma^*}\; .
$$

Therefore one has
$$
\ww(g_1) - \ww(G^i_2)-\ww(G^j_1) = \ww(G^i_1/G^i_2) +
\ww(G^0_1)
$$
and this point is strictly under the line $\cal L$ in
$\R^2\supset\Z_{\ge 0}^2$ which goes through the origin and the
point
$(m_i^{\sigma^*}, m_j^{\sigma^*})$.

On the other hand
$$
\ww(g_1) - \ww(G^i_2)-\ww(G^j_1) =
\ww(g_2) - \ww(G^i_2)-\ww(G^j_2) = \ww(G^0_2).
$$
The last point lies on the line $\cal L$. This proves the
statement.
\end{proof}

\begin{lemma}\label{lemma6}
Let $\kk_1, \ldots , \kk_p$ be different multi-indices
such that there exist $y_i\in Y_{\kk_i}$, for all $i=1,\ldots ,
p$ with $\w{\Pi}(y_1) = \cdots=\w{\Pi}(y_p)$.
Then there exist maximal cuts $\sigma_i$ of
$\kk_i$, $i=1,\ldots, p$, which are comparable with each other,
i.e.,
all of them lie on one and the same geodesic from the vertex ${\bf 1}$
to a strict transform $\widetilde C_j$ of a branch of the curve $C$.
\end{lemma}

\begin{proof}
First let us prove the existence of cuts of all the multi-indices
$\kk_i$, $i=1\ldots, p$. It is sufficient to prove this for $p=2$.
By Lemma~{\ref{lemma3}}, there exists a cut, say of $\kk_1$, at a
vertex $\sigma_1\in \Delta^\prime$. If there  is no cut of
$\kk_2$, then Lemma~{\ref{lemma5}} says that there exists a cut
$\sigma_i$ of $\kk_1$ on the geodesic from $st_{1}$ to $\w{C_i}$
for all $i=1,\ldots , r$. Using Lemma~{\ref{lemma4}} and the
Remark~1 after it, one gets that $$ \vv(g_1/g_2; \ge st_1)\ge
\mm^{st_1} +\w{\mm}^{st_1}\; . $$ One has (see, e.g., \cite{D1} or
item (2) in \ref{irred})
$$ \w{m}^{st_1}\ge
\mm^{\beta_q} = (n_q+1)\cdot\mm^{\alpha_q}
 > n_q\cdot\mm^{\alpha_q} +(n_{q-1}+1)\cdot\mm^{\alpha_{q-1}}
 >\cdots
$$
$$\cdots >
\sum\limits_{p=1}^q n_p \cdot\mm^{\alpha_p} + \mm^{\alpha_0}
 >\sum\limits_{p=1}^q n_p \cdot\mm^{\alpha_p} - \mm^{\alpha_0}\;
.
$$
Since $\mm=\mm^{st_1}+\sum\limits_{p=1}^q n_p \cdot\mm^{\alpha_p} -
\mm^{\alpha_0}$ is the maximal element of (the finite set) $S^\prime_1$,
$\mm^{st_1}+\w{\mm}^{st_1}$ is strictly bigger than $\mm$.
Therefore
$$
\vv(g_1)-\vv(g_2)=\vv(g_1/g_2;\ge {\bf 1}) > 0
$$
what contradicts the supposition that $\vv(g_1)=\vv(g_2)$.

To prove the existence of maximal cuts which are comparable with each
other, we shall rather prove the following statement. Let
$\sigma\in\Delta^\prime$ and suppose that among the vertices $\sigma^\prime
\ge \sigma$, $\sigma^\prime\in\Delta^\prime$, there are cuts of all the
multi-indices $\kk_i$, $i=1, 2, \ldots, p$. Then there exist
maximal cuts
$\sigma_i$ of $\kk_i$, $\sigma_i\ge \sigma$, which are comparable with
each other. We shall use the simultaneous induction both on the number
of multi-indices $p$ and on the vertex $\sigma$ (in the inverse
order).
If $p=1$ or if $\sigma$ is a maximal element in $\Delta^\prime$,
the statement
is trivial (in the last case all cuts coincide with $\sigma$). Let $Q_0$,
$Q_1$, \dots, $Q_{s-1}$ ($s=s_\sigma$) be essential points of the component
$E_\sigma$; if $\sigma\ne st_1$, we assume that the component of
$D^\prime\setminus\widetilde E_\sigma$ corresponding to the essential point
$Q_0$ contains the divisor $E_{st_1}$. Let $\Gamma(\sigma,\ell)$ be the
subgraph of $\Gamma$ corresponding to the essential point $Q_\ell$, $\ell=0,
1, \dots, s-1$. One can meet one of the three following situations. \newline
1) $\sigma$ is a maximal cut of one of the multi-indices $\kk_i$, say, of
$\kk_p$. In this case the statement follows from one for
multi-indices $\kk_1$,
\dots, $\kk_{p-1}$. \newline
2) $\sigma$ is not a maximal cut of any of the multi-indices $\kk_i$, and there
exists $\ell$, $1\le \ell\le s-1$ ($0\le \ell\le s-1$ if $\sigma=st_1$),
such that the subgraph $\Gamma(\sigma, \ell)$ contains cuts of all the
multi-indices $\kk_i$. In this case the statement follows from one applied to
the first (i.e., the minimal) vertex in $\Gamma(\sigma,\ell)$. \newline
3) $\sigma$ is not a maximal cut of any of the multi-indices $\kk_i$, and there
is no $\ell$, $1\le \ell\le s-1$ ($0\le \ell\le s-1$ if $\sigma=st_1$),
such that the subgraph $\Gamma(\sigma, \ell)$ contains cuts of all the
multi-indices $\kk_i$. In this case there exist $\ell_1$, $\ell_2$
($1\le \ell_u\le s-1$; $0\le \ell_u\le s-1$ if $\sigma=st_1$), $i_1$, and
$i_2$ such that $\Gamma(\sigma, \ell_1)$ contains a cut of $\kk_{i_1}$, but
does not contain a cut of $\kk_{i_2}$, and vice versa $\Gamma(\sigma, \ell_2)$
contains a cut of $\kk_{i_2}$, but does not contain a cut of $\kk_{i_1}$.
Without loss of generality one can suppose that $i_1=1$, $i_2=2$. Let $i$
(respectively $j$) be such that the strict transform $\widetilde C_i$ lies in
$\Gamma(\sigma,\ell_1)$ (respectively $\widetilde C_j$ lies in
$\Gamma(\sigma,\ell_2)$) and the geodesic from $\sigma$ to $\widetilde C_i$
contains a maximal cut $\sigma_1$ of $\kk_1$ (respectively the geodesic from
$\sigma$ to $\widetilde C_j$ contains a maximal cut $\sigma_2$ of $\kk_2$).

Now we shall use the same arguments (and the same notations) as in
the proof of Lemma~{\ref{lemma5}} applied to the vertex $\sigma$,
$\Gamma(\sigma,\ell_1)$ and $\Gamma(\sigma,\ell_2)$.

Let $g_{u} = G^0_u\cdot G^i_u \cdot G^j_u$ where
$G^i_u = \prod\limits_{\tau\in \Gamma(\sigma, \ell_1)}
\varphi^{(\tau)}_u$,
$G^j_u = \prod\limits_{\tau\in \Gamma(\sigma, \ell_2)}
\varphi^{(\tau)}_u$,
$G^0_u = g_u/(G^i_u\cdot G^j_u)$, $u=1, 2$.

Since $\sigma_1$ (respectively $\sigma_2$) is a cut of $\kk_1$
(respectively
of $\kk_2$) and $\Gamma(\sigma,\ell_1)$ (respectively
$\Gamma(\sigma,\ell_2)$) does not contains cuts of
$\kk_2$ (respectively of
$\kk_1$), using repeatedly Lemma~{\ref{lemma4}}
and the Remark after it one can show that
$$
\begin{array}{rl}
\ww(G^i_1/G^i_2) &\ge (\widetilde m_i^{\sigma},
\widetilde m_j^{\sigma}) > (0,0)\,, \\
(v_j(G^i_1)-v_j(G^i_2))m_i^{\sigma} &<
(v_i(G^i_1)-v_i(G^i_2))m_j^{\sigma}\,, \\
\ww(G^j_2/G^j_1) &\ge (\widetilde m_i^{\sigma},
\widetilde m_j^{\sigma}) > (0,0)\,, \\
(v_j(G^j_2)-v_j(G^j_1))m_i^{\sigma} &>
(v_i(G^j_2)-v_i(G^j_1))m_j^{\sigma}\,. \\
\end{array}
$$
One has
$$
\ww(g_1) - \ww(G^i_2)-\ww(G^j_1) = \ww(G^i_1/G^i_2) +
\ww(G^0_1)
$$
and this point is strictly under the line $\cal L$ in
$\R^2\supset\Z_{\ge 0}^*$ which goes through the origin
and the point $(m_i^{\sigma}, m_j^{\sigma})$.
On the other hand,
$$
\ww(g_1) - \ww(G^i_2)-\ww(G^j_1) =
\ww(g_2) - \ww(G^i_2)-\ww(G^j_1) = \ww(G^j_2/G^j_1) +
\ww(G^0_2)
$$
where this point is strictly over the line $\cal L$.
This proves the statement.
\end{proof}

\begin{remarks}
{\bf 1.} Under  the conditions of the Lemma~{\ref{lemma6}}, let
$\sigma_1$, \dots,
$\sigma_p$ be maximal cuts of multi-indices $\kk_1$, \dots,
$\kk_p$ such that
$\sigma_1\le \sigma_2\le \cdots \le\sigma_p$. Then, by
Lemma~{\ref{lemma5}},
on each geodesic from the vertex $\sigma_1$ to a strict transform $\w{C}_j$
with $\w{C}_j>\sigma_1$ there exists a cut of $\kk_i$ for $i=2, \ldots, p$.

\noindent {\bf 2.} As a consequence of Lemma~{\ref{lemma6}}, if
a multi-index $\kk$ has no cuts then the image $\w{\Pi}(Y_{\kk})$ of the
corresponding component $Y_{\kk}$ of the space $\widetilde Y$ does not
intersects the image $\w{\Pi}(Y_{\kk'})$ for $\kk'\neq \kk$.
\end{remarks}

\begin{lemma}\label{lemma7}
Let $Q_0$, $Q_1$, \dots, $Q_{s-1}$ be (different) points of a
projective line $E$, $\widetilde{E} = E- \{Q_\ell : \ell=0, 1,
\ldots, s-1\}$, let $P^o_1$, \dots, $P^o_k$ be $k$ points (not
necessarily different), different from $Q_0$, $Q_1$, \dots,
$Q_{s-1}$. Let $\Phi$ be the map from
$S^k \widetilde{E}$ to $\P ((\C^*)^s)$ defined in the following way.
For an element from $S^k\widetilde{E}$; i.e., for $k$ points
$P_1$, \dots, $P_k$, let $\psi$ be a meromorphic function on $E$ with
zeroes at the points $P_1$, \dots, $P_k$ and poles at the points
$P^o_1$, \dots, $P^o_k$;
let $\Phi(\{P_j\}) := (\psi(Q_0): \psi(Q_1): \ldots :
\psi(Q_{s-1}))$.
Then, if $k\ge s-1$, one has $\mbox{Im }\Phi = \P((\C^*)^s)$; if
$k\le s-1$, $\Phi$ is an embedding. Moreover in both cases $\Phi$
is a (smooth) locally trivial (in fact a trivial) fibration over its image the
fibre of which is a (complex) affine space of dimension $\max(0, k-s+1)$.
\end{lemma}

\begin{proof}
Without loss of generality one can suppose that
$P^o_1=P^o_2=\cdots = P^o_k=P^o$. Let us choose an affine
coordinate on $E$ such that $P^o=\infty$. Then $\psi$ is a
polynomial of degree $\le k$ with zeroes at those points
$P_1$, \dots, $P_k$ which are different from $P^o$. Let $z_\ell$ be the
coordinate of the point $Q_\ell$, $\ell=0, 1, \ldots, s-1$.

For $k\ge s-1$, the statement that the map $\Phi$ is onto can be reduced
to the following obvious one: for an arbitrary prescribed set of values
$\{\psi_0$, $\psi_1$, \dots, $\psi_{s-1}\}$, there exists a
polynomial $\psi$ of degree $\le k$ such that
$\psi(z_\ell)=\psi_\ell$, $\ell=0, 1, \ldots, s-1$. The statement that
$\Phi$ is a locally trivial fibration over its image
follows from the fact that if $\psi_1$ and $\psi_2$ are polynomials
with coinciding values at the points $Q_\ell$, $\ell=0, 1, \ldots, s-1$,
then $\psi_1=\psi_2+q(z)(z-z_0)(z-z_1)\cdot\ldots\cdot(z-z_{s-1})$
where $q(z)$ is an arbitrary polynomial of degree $k-s$.

For $k\le s-1$, the statement follows from the fact that such a
polynomial of degree $\le s-1$ is unique.
\end{proof}

\subsection{Free action of the torus.}\label{s-torus}
Let $\sigma\in\Delta^\prime$, $s=s_\sigma$. There is defined the
following free action $T=T(\sigma)$ of the group
$\P((\C^*)^s)\cong(\C^*)^{s-1}$ on
$\Z_{\ge 0}^r\times\P((\C^*)^r)$. Let
$Q_0$,
$Q_1$, \dots $Q_{s-1}$ be essential points on the component
$E_\sigma$ and, for $1\le i\le r$, let $Q_{\ell(i)}$ be the
essential point corresponding to the connected component of
the complement $(f\circ\pi)^{-1}(0)\setminus {\stackrel{\circ}{E}}_\sigma$
which contains the strict transform $\widetilde C_i$. Let
$\underline{c}=(c_0:c_1:\ldots:c_{s-1})\in \P((\C^*)^s)$. Then,
for $(\vv, \aa)=(v_1,\ldots,v_r; a_1:\ldots:a_r)\in \Z_{\ge
0}^r\times\P((\C^*)^r)$,
one has $T_{\underline c}(\vv, \aa):=(\vv;
c_{\ell(1)}\cdot a_1:\ldots:c_{\ell(r)}\cdot a_r)$ (i.e., all
coordinates $a_i$ of $\aa$ such that the strict transform
$\widetilde C_i$ intersects the component of $D^\prime\setminus
{\stackrel{\circ}{E}}_\sigma$ corresponding to one essential point $Q_\ell$ are
multiplied by one and the same number $c_\ell$).

\begin{corollary}\label{cor2}
Suppose $\sigma$ is a cut of $\kk$ and $y\in Y_\kk$,
$y=[K_\sigma]\cdot\tt^{k_\sigma\mm^\sigma}\times y^*$
($\#K_\sigma=k_\sigma$). Then $\widetilde\Pi(S^{k_\sigma}
\widetilde E_\sigma\cdot\tt^{k_\sigma\mm^\sigma}\times y^*)$
is the orbit of $\widetilde\Pi(y)$ under the described action
$($and thus is homeomorphic to $\P((\C^*)^s)\cong(\C^*)^{s-1}$$)$.
\end{corollary}

\begin{proof}
Indeed, let $K_\sigma^\prime$ be a subset of $\widetilde E_\sigma$
with $k_\sigma$ elements (i.e., $[K_\sigma^\prime]$ is an
element of $S^{k_\sigma}\widetilde E_\sigma$), let
$y^\prime=[K_\sigma^\prime]\cdot\tt^{k_\sigma\mm^\sigma}\times y^*$,
let $g$ and $g^\prime$ be functions from ${\cal O}_{\C^2, 0}$
corresponding to the points $y$ and $y^\prime$ of $Y_\kk$, and
let $\psi=\widetilde g^\prime/\widetilde g$, where
$\widetilde g=g\circ\pi$ and $\widetilde g^\prime=g^\prime\circ\pi$
are the liftings of the functions $g$ and $g^\prime$ to the space
$X$ of the resolution. Then $\psi_{\vert E_\sigma}$ is a
meromorphic function on the projective line $E_\sigma$ with
$k_\sigma$ zeroes at the points of the set $K_\sigma^\prime$ and
$k_\sigma$ poles at the points of the set $K_\sigma$ (such a
function is well-defined up to the multiplication by a constant).
Moreover $\psi$ is constant on each connected component of
$D^\prime\setminus \widetilde E_\sigma$ and its value on this component
coincides with the value of $\psi$ at the corresponding essential
point of $E_\sigma$. Therefore the statement follows from
Lemma~{\ref{lemma7}} (for the case $k\ge s-1$).
\end{proof}

\begin{corollary}\label{cor3}
Suppose that $y\in Y_\kk$ and $\kk$ has a cut on each geodesic
from the vertex
$\sigma$ to a strict transform $\widetilde C_i$ with $\widetilde C_i > \sigma$.
Then $\widetilde\Pi(Y_\kk)$ contains the orbit of
$\widetilde\Pi(y)$ under the described action.
\end{corollary}

\begin{proof}
Let $\sigma_1$, \dots, $\sigma_p$ be the minimal elements in the
set of cuts of $\kk$ which are $\ge\sigma$. If $\sigma_1=\sigma$
(and thus $p=1$), the statement is obvious (see
Corollary~{\ref{cor2}}). Let $y=\prod\limits_{j=1}^p
([K_{\sigma_j}]\cdot\tt^{k_{\sigma_j}\mm^{\sigma_j}})\times y^*$
($K_{\sigma_j}\subset\widetilde E_{\sigma_j}$, $\# K_{\sigma_j}=k_{\sigma_j}$).
For $j=1, \ldots, p$, let $Q_0^j$, $Q_1^j$, \dots, $Q_{s_j-1}^j$
($s_j=s_{\sigma_j}$) be the essential points of the component
$E_{\sigma_j}$ numbered so that the connected component of
$D^\prime\setminus \widetilde E_{\sigma_j}$, corresponding to the point
$Q_0^j$, contains the component $E_\sigma$ of the exceptional
divisor (or equivalently the first separation component
$E_{st_1}$). Let $\widetilde\ell(j)$ be such that the connected
component of $D^\prime\setminus\widetilde E_\sigma$ corresponding to
the essential point $Q_{\widetilde\ell(j)}\in\widetilde E_\sigma$
contains the component $E_{\sigma_j}$ of the essential divisor.
One knows (Lemma~{\ref{lemma7}}) that, for any set
of non-zero numbers $q_0^j$, $q_1^j$, \dots, $q_{s_j-1}^j$, there exists a
subset
$K_{\sigma_j}^\prime\subset\widetilde E_{\sigma_j}$ with
$\# K_{\sigma_j}^\prime=k_{\sigma_j}$ such that a meromorphic
function $\psi_j$ on $E_{\sigma_j}$ with zeroes at the points of
the set $K_{\sigma_j}^\prime$ and poles at the points of
the set $K_{\sigma_j}$ has values $q_0^j$, $q_1^j$, \dots, $q_{s_j-1}^j$
at the points $Q_0^j$, $Q_1^j$, \dots, $Q_{s_j-1}^j$
respectively. For $\underline c=(c_0 : c_1 : \ldots : c_{s-1})
\in \P((\C^*)^s)$,
let $K_{\sigma_j}^\prime$ be such that $\psi_j(Q_0^j)=1$,
$\psi_j(Q_m^j)=c_{\widetilde\ell(j)}/c_0$ for $m\ge 1$, and let
$y^\prime=\prod\limits_{j=1}^p ([K_{\sigma_j}^\prime]\cdot
\tt^{k_{\sigma_j}\mm^{\sigma_j}})\times y^*$. Then
$\widetilde\Pi(y^\prime)=T_{\underline c}(\widetilde\Pi(y))$.
\end{proof}

\begin{lemma}\label{lemma8}
For a multi-index $\kk=\{k_\delta,\, \mm,\, k_\sigma\}$, the
map $\widetilde\Pi_{\vert Y_{\kk}}: Y_{\kk}\to \widetilde\Pi(Y_{\kk})$
is a locally trivial fibration.
$($Note that $\widetilde\Pi(Y_\kk)\subset \P F_\vv$, where $\vv=\vv(\kk)$$)$.
\end{lemma}

\begin{proof}
Let
$$
y = \prod\limits_{\sigma\in\Delta^\prime}
\left([K_{\sigma}^o]\cdot\tt^{k_{\sigma}\mm^\sigma}\right)\times
(\bullet\cdot\tt^{\mm})\times\prod\limits_{\delta\in \Delta}
(\bullet\cdot\tt^{k_{\delta}\mm^\delta})
$$
be a point of $Y_\kk$. For $\sigma\in\Delta^\prime$ let
$\Phi_\sigma:S^{k_\sigma}{\widetilde{E}}_\sigma\to\P((\C^*)^{s_\sigma})$
be the map ($\Phi$) described in Lemma~\ref{lemma7} for $\widetilde E =
\widetilde E_\sigma$, $k=k_\sigma$, $\{P_1^o, \ldots, P_k^o\}=K_{\sigma}^o$,
and $\{Q_0, Q_1, \ldots, Q_{s-1}\}=
\{Q_0^\sigma,Q_1^\sigma, \ldots, Q_{s_\sigma-1}^\sigma\}$, let
$\phi_\sigma=\Phi_\sigma\circ\pi_\sigma$ where $\pi_\sigma:Y_\kk\to
S^{k_\sigma}{\widetilde{E}}_\sigma$ is the natural projection, and let
$\Psi=\prod_{\sigma}\phi_\sigma:Y_\kk\to \prod_\sigma\P((\C^*)^{s_\sigma})$.
 From Lemma~\ref{lemma7} it follows that $\Psi$ is a locally trivial fibration
over its image. Let $M: \prod_\sigma\P((\C^*)^{s_\sigma}) \to \P F_\vv$ be the
map defined in the following way.
For $1\le i\le r$, $\sigma\in\Delta^\prime$, let $Q_{\ell_\sigma(i)}^\sigma$ be
the esential point of the component $E_\sigma$ of the exceptional divisor,
corresponding to the connected component of the complement
$(f\circ\pi)^{-1}(0)\setminus {\stackrel{\circ}{E}}_\sigma$
which contains the strict transform $\widetilde C_i$.
For $\underline c^\sigma=(c_0^\sigma : c_1^\sigma : \ldots:
c_{s_\sigma-1}^\sigma)
\in\P((\C^*)^{s_\sigma})$, $\sigma\in\Delta^\prime$,
$M(\prod_\sigma \underline c^\sigma):=\left(\prod_\sigma
c_{\ell_\sigma(1)}^\sigma : \ldots : \prod_\sigma c_{\ell_\sigma(r)}^\sigma,
\right)\cdot\widetilde\Pi(y)\in\P F_{\vv}$. The map $M$ is a locally
trivial fibration. Now the statement follows from the fact that
$\widetilde\Pi_{\vert Y_{\kk}}=M\circ\Psi$.
\end{proof}

\begin{statement}
For all multi-indices $\kk$ with $\vv(\kk)=\vv$ (there is
a finite number of them), one can
construct subspaces $Y_{\kk}^\prime\subset Y_{\kk}$ such that \newline
1) $\widetilde\Pi(\bigcup Y_{\kk}^\prime)=\widetilde\Pi(\bigcup Y_{\kk})$;
\newline
2) $\widetilde\Pi(Y_{\kk_1}^\prime)$ does not intersect
$\widetilde\Pi(Y_{\kk_2}^\prime)$ for $\kk_1\ne\kk_2$; \newline
3) $\chi(Y_{\kk}\setminus Y_{\kk}^\prime)=0$; \newline
4) either $\widetilde\Pi_{\vert Y_{\kk}}$ is one--to--one on its
image, or $\chi(Y_{\kk}^\prime)=0$ and $\chi(\widetilde\Pi(Y_{\kk}^\prime))=0$.
\end{statement}

\begin{proof}
Let us order multi-indices $\kk$ with $\vv(\kk)=\vv$ in an
arbitrary way. For $\kk_1$ with $\vv(\kk_1)=\vv$, let ${\cal
I}_{\kk_1}=\{(\vv,\aa)\in\P F_\vv: \exists y_1\in Y_{\kk_1},\,
\exists\kk_2>\kk_1,\,\exists y_2\in
Y_{\kk_2}:\,\widetilde\Pi(y_1)= \widetilde\Pi(y_2)= (\vv,\aa)\}\;
$,  $Z_{\kk_1}=\widetilde\Pi^{-1}({\cal I}_{\kk_1})\cap
Y_{\kk_1}$, $Y_{\kk_1}^\prime=Y_{\kk_1}\setminus Z_{\kk_1}$. The
Euler characteristic of the subspace ${\cal I}_{\kk_1}$ is an
alternative sum of Euler characteristics of the subspaces ${\cal
I}(\kk_1, \kk_2, \ldots, \kk_p)=\{(\vv,\aa)\in\P F_\vv: \exists
y_i\in Y_{\kk_i},\, i=1, 2, \ldots, p:\,\widetilde\Pi(y_1)=
\widetilde\Pi(y_2)=\ldots=\widetilde\Pi(y_p)=(\vv,\aa)\}$ with
$p\ge 2$.

If the set ${\cal I}(\kk_1, \kk_2, \ldots, \kk_p)$ is not empty,
according to Lemma~{\ref{lemma6}} there exist maximal cuts
$\sigma_i$ of $\kk_i$, $i=1, 2, \ldots, p$, which are comparable
with each other, i.e., which lie on one and the same geodesic in
the graph $\Gamma$ from $st_1$ to a strict transform of a branch
of the curve $C$. Let $\sigma_{i_0}$ be (one of) the smallest of
these cuts (all smallest cuts coincide with each other). The
remark after Lemma~{\ref{lemma6}} says that for each $i\ne i_0$ on
each geodesic from the vertex $\sigma_{i_0}$ to a strict transform
$\widetilde C_j$ with $\widetilde C_j > \sigma_{i_0}$ there exists
a cut of $\kk_i$. By Corollary~{\ref{cor3}} the (semianalytic)
subspace ${\cal I}(\kk_1, \kk_2, \ldots, \kk_p)$ is invariant
with respect to the described above free action of the group
$\P((\C^*)^s)\cong(\C^*)^{s-1}$, where $s=s_{\sigma_{i_0}}$.
Therefore $\chi({\cal I}(\kk_1, \kk_2, \ldots, \kk_p))=0$,
$\chi({\cal I}_{\kk_1})=0$, and (since the map
$\widetilde\Pi_{\vert Y_{\kk_1}}$ is a locally trivial fibration)
$\chi(Z_{\kk_1})=0$.

Obviously the sets $Y_\kk^\prime$ satisfy the conditions 1)~-- 3).
If $\widetilde\Pi_{\vert Y_{\kk}}$ is not one--to--one on its
image, there exists a cut $\sigma$ of $\kk$. In this case the
space $Y_\kk$ is a product of a space and the symmetric power
$S^{k_\sigma}\widetilde E_\sigma$ with $k_\sigma\ge s_\sigma-1$.
Therefore $\chi(Y_\kk)=0$ and $\chi(Y_\kk^\prime)=0$. The image
$\widetilde\Pi(Y_\kk)$ is invariant with respect to the free
action of the group
$\P((\C^*)^{s_\sigma})\cong(\C^*)^{s_\sigma-1}$.
Therefore $\chi(\widetilde\Pi(Y_\kk))=0$,
$\chi(\widetilde\Pi(Y_\kk^\prime))= \chi(\widetilde\Pi(Y_\kk))-\chi({\cal
I}_\kk)=0$. \end{proof}

The Statement obviously implies Proposition~\ref{prop4} and thus
Theorem~\ref{theo4} has been proved.


\begin{thebibliography}{12}
\bibitem[11]{Eb} Ebeling W. Poincar\'e series and monodromy of
a two--dimensional quasihomogeneous hypersurface singularity.
To be published in manuscripta math.
\end{thebibliography}
\end{document}